\documentclass[reqno,11pt,a4paper]{amsart}
\usepackage{amsthm}
\usepackage{amssymb}
\usepackage{cite}
\usepackage[all]{xy}
\usepackage[full]{textcomp}
\usepackage{newpxtext}
\usepackage{cabin}
\usepackage[varqu,varl]{inconsolata}
\usepackage[bigdelims,vvarbb]{newpxmath}
\usepackage[cal=cm,bb=esstix,bbscaled=1.05]{mathalfa}
\usepackage[margin=2.7cm]{geometry}
\usepackage{graphicx}
\usepackage[usenames,dvipsnames]{color}
\usepackage{hyperref}
\hypersetup{
    colorlinks,
    linkcolor={red!50!black},
    citecolor={blue!50!black},
    urlcolor={blue!80!black}
}
\usepackage{longtable} 
\usepackage{colortbl} 
\usepackage{arydshln} 
\usepackage{booktabs} 
\usepackage{mleftright}
\usepackage{appendix}
\usepackage{tikz}
\usetikzlibrary{positioning}
\usepackage{tikz-3dplot}
\usepackage{cite}
\usepackage{pdflscape}
\usepackage{bm}
\usepackage{comment}

\AtBeginDocument{%
\def\MR#1{}
}
\theoremstyle{plain} 
\newtheorem{Lemma}{Lemma}[section]
\newtheorem{Thm}[Lemma]{Theorem}
\newtheorem{Prop}[Lemma]{Proposition}
\newtheorem{Cor}[Lemma]{Corollary}
\newtheorem{Conj}[Lemma]{Conjecture}
\newtheorem{Que}{Question}

\theoremstyle{definition}
\newtheorem{Defn}[Lemma]{Definition}
\newtheorem{Ex}[Lemma]{Example}
\newtheorem{Convention}[Lemma]{Convention}

\theoremstyle{remark}
\newtheorem{Rem}[Lemma]{Remark}

\numberwithin{equation}{section}

\newcommand{\ZZ}{\mathbb{Z}}
\newcommand{\PP}{\mathbb{P}}
\newcommand{\CC}{\mathbb{C}}
\newcommand{\QQ}{\mathbb{Q}}

\newcommand{\FF}{\mathbb{F}}

\newcommand{\Cstar}{\CC^\times}

\newcommand{\calG}{\mathcal{G}}

\newcommand{\cR}{\mathcal{R}}

\newcommand{\GIT}{/\!\!/}

\let\ker\undefined

\DeclareMathOperator{\ker}{ker}

\renewcommand{\emptyset}{\varnothing}

\newcommand{\bx}{\bm{x}}

\newcommand{\NQ}{N_{\QQ}}
\newcommand{\MQ}{M_{\QQ}}
\newcommand{\orig}{\textbf{0}}
\newcommand{\abs}[1]{\left\vert{#1}\right\vert}
\DeclareMathOperator{\GL}{GL}
\DeclareMathOperator{\SL}{SL}

\DeclareMathOperator{\Hom}{Hom}

\DeclareMathOperator{\coeff}{coeff}
\DeclareMathOperator{\linspan}{span}
\DeclareMathOperator{\id}{Id}
\newcommand{\bdry}[1]{\partial{#1}}
\newcommand{\intr}[1]{{#1}^\circ}
\newcommand{\Vol}[1]{\operatorname{Vol}\mleft({#1}\mright)}
\newcommand{\V}[1]{\operatorname{vert}\mleft({#1}\mright)}
\renewcommand{\dim}[1]{\operatorname{dim}\mleft({#1}\mright)}

\renewcommand{\gcd}[1]{\operatorname{gcd}\mleft\{{#1}\mright\}}

\newcommand{\Newt}[1]{\operatorname{Newt}\mleft({#1}\mright)}
\newcommand{\conv}[1]{\operatorname{conv}\mleft({#1}\mright)}
\newcommand{\sconv}[1]{\operatorname{conv}\mleft\{{#1}\mright\}}

\newcommand{\cB}{\mathcal{B}}
\newcommand{\SC}[1]{\operatorname{SC}\mleft({#1}\mright)}
\newcommand{\res}[1]{\operatorname{res}\mleft({#1}\mright)}
\newcommand{\MM}[2]{\mathrm{MM}_{#1\text{--}#2}} 
\newcommand{\evnrow}{\rowcolor[gray]{0.95}}
\newcommand{\oddrow}{}
\newcommand{\NotMinkowski}[1]{\color{blue}#1}

\graphicspath{{images/}}
\begin{document}
\author[T. Coates]{Tom Coates}
\address{Department of Mathematics\\Imperial College London\\180 Queen's Gate\\London\\SW7 2AZ\\UK}
\email{t.coates@imperial.ac.uk}
\author[A.\,M.\,Kasprzyk]{Alexander M.~Kasprzyk}
\address{School of Mathematical Sciences\\University of Nottingham\\Nottingham\\NG7 2RD\\UK}
\email{a.m.kasprzyk@nottingham.ac.uk}
\author[G.\,Pitton]{Giuseppe Pitton}
\address{Department of Mathematics\\Imperial College London\\180 Queen's Gate\\London\\SW7 2AZ\\UK}
\email{g.pitton@imperial.ac.uk}
\author[K.\,Tveiten]{Ketil Tveiten}
\address{Department of Mathematics\\Uppsala University\\Box 256 SE-$751$\ $05$\ Uppsala\\Sweden}
\email{ketiltveiten@gmail.com}
\title{Maximally Mutable Laurent Polynomials}
\begin{abstract}
    We introduce a class of Laurent polynomials, called maximally mutable Laurent polynomials (MMLPs), that we believe correspond under mirror symmetry to Fano varieties. A subclass of these, called rigid, are expected to correspond to Fano varieties with terminal locally toric singularities. We prove that there are exactly 10 mutation classes of rigid MMLPs in two variables; under mirror symmetry these correspond one-to-one with the 10 deformation classes of smooth del~Pezzo surfaces. Furthermore we give a computer-assisted classification of rigid MMLPs in three variables with reflexive Newton polytope; under mirror symmetry these correspond one-to-one with the 98 deformation classes of three-dimensional Fano manifolds with very ample anticanonical bundle. We compare our proposal to previous approaches to constructing mirrors to Fano varieties, and explain why mirror symmetry in higher dimensions necessarily involves varieties with terminal singularities. Every known mirror to a Fano manifold, of any dimension, is a rigid MMLP. 
\end{abstract}
\maketitle
\section{Introduction}
Recently a new approach to the classification of Fano manifolds was proposed~\cite{ProcECM}. This centres on mirror symmetry, in the form of a conjectural relationship between Fano manifolds and Laurent polynomials. An $n$-dimensional Fano manifold~$X$ determines the \emph{regularised quantum period}
\begin{equation}
    \label{eq:regularised quantum period}
    \widehat{G}_X(t)=1 + \sum_{k=2}^\infty k!c_kt^k.
\end{equation}
This is a generating function for certain Gromov--Witten invariants of $X$: roughly speaking, the coefficient $c_k$ is the number of degree-$k$ rational curves on $X$ that pass through a given point (see~\cite{ProcECM} for details). The quantum period is expected to characterise $X$. Mirror symmetry suggests that the power series $\widehat{G}_X$ also arises from a Laurent polynomial $f\in\CC[x_1^{\pm1},\ldots,x_n^{\pm1}]$, and when this happens there is a close connection between the geometry of $X$ and of $f$. More precisely, given a Laurent polynomial $f\in\CC[x_1^{\pm 1},\ldots,x_n^{\pm 1}]$ one can form its \emph{classical period}:
\begin{equation}\label{eq:period_integral}
\pi_f(t)=\left(\frac{1}{2\pi i}\right)^n\int_{\abs{x_1}=\ldots=\abs{x_n}=1}\frac{1}{1-tf}\frac{dx_1}{x_1}\cdots\frac{dx_n}{x_n},\qquad t\in\CC,\abs{t}\ll\infty.
\end{equation}
The Taylor expansion of $\pi_f$, which we also denote $\pi_f$, can be readily calculated~\cite{ACGK12,Lairez}:
\begin{equation}\label{eq:period_sequence}
\pi_f(t)=\sum_{k=0}^\infty\coeff_1(f^k)t^k.
\end{equation}
Here $\coeff_1(f^k)$ denotes the coefficient of the constant term of $f^k$. If
\[
\widehat{G}_X=\pi_f
\]
then we say that $f$ is a \emph{mirror partner} for $X$. This is a very strong condition on $f$. However certain Laurent polynomials are known to arise as mirror partners to $n$-dimensional Fano manifolds~\cite{ProcECM,QC105,CGKS14}, and in particular mirror partners exist for all Fano manifolds of dimension up to three.

Henceforth let us consider only Laurent polynomials $f$ such that the exponents of monomials in $f$ generate the lattice $\ZZ^n$; this will avoid a subtlety about finite coverings of toric varieties. If a Laurent polynomial $f$ is a mirror partner to a Fano manifold $X$ then it is expected that $X$ admits a \emph{$\QQ$-Gorenstein (qG-) degeneration} to the singular toric variety $X_P$ associated to the Newton polytope $P=\Newt{f}$ of $f$~\cite{AK14,pragmatic}. Here $X_P$ is given by taking the \emph{spanning fan}, also called the face fan or central fan, whose cones span the faces of~$P$. Given a mirror partner $f$, new mirror partners can be generated via \emph{mutation}; this process, which produces new Laurent polynomials with the same classical period, is described in detail below. Suppose that $f$ is a mirror partner to $X$. It is conjectured that the Laurent polynomial $g$ is also a mirror partner to $X$ if and only if $f$ and $g$ are connected via a sequence of mutations. Such a connection implies in particular that the corresponding toric varieties $X_P$ and $X_Q$, $P=\Newt{f}$, $Q=\Newt{g}$, are related via qG-deformation~\cite{Ilten}.

\subsection{Fano polytopes}
We now attempt to reverse the construction described above. Let $N$ be a lattice of rank $n$, and let $P\subset\NQ:=N\otimes_\ZZ\QQ$ be a convex lattice polytope.

\begin{Que}\label{qu:assign_laurent}
Does there exist a Laurent polynomial
\[
f=\sum_{v\in P\cap N}c_v\bx^v,
\]
with coefficients $c_v$ such that $\Newt{f}=P$, with $f$ a mirror partner to some Fano manifold?
\end{Que}

Without loss of generality, we may make the following three assumptions.
\begin{enumerate}
\item\label{item:origin_in_intr}
We may assume that the origin $\orig$ of $N$ is contained in the relative interior $\intr{P}:=P\setminus\bdry{P}$ of $P$: if $\orig\notin P$ then, for any Laurent polynomial $f$ with $\Newt{f}=P$, the period $\pi_f$ will be constant, and hence cannot be the quantum period of a Fano manifold; if $\orig\in\bdry{P}$ then we can reduce to a lower-dimensional situation by considering the smallest-dimensional face $F$ of $P$ with $\orig\in F$.
\item
We can assume that $\dim{P}=n$, since otherwise we can restrict to the sublattice given by $\linspan_\QQ(P)\cap N$.
\item\label{item:generate_lattice}
Consider the sublattice
\[
N':=\sum_{v\in P\cap N}v\cdot\ZZ
\]
of $N$ generated by the lattice points of $P$. Restricting to the sublattice $N'$ will not change the corresponding period $\pi_f$, since $\coeff_1(f^k)$ depends on $N$ only via the linear relations between the lattice points of $P$. Hence we may assume that $N'=N$.
\end{enumerate}

Recall that we want to associate to $P$ a toric variety $X_P$ via the spanning fan. The rays of this fan are spanned by the vertices $\V{P}$ of $P$, so we require that the vertices are primitive lattice points of $N$ (and hence correspond to the lattice generators of the rays). By assumption~\eqref{item:origin_in_intr} the fan is complete, and we see that $X_P$ is an $n$-dimensional toric Fano variety.

\begin{Defn}
A convex lattice polytope $P\subset\NQ$ is called \emph{Fano} if it is of maximum dimension with respect to the ambient lattice $N$, if it contains the origin $\orig$ in its interior $\intr{P}$, and if the vertices of $P$ are primitive lattice elements. 
\end{Defn}

A Fano polytope $P$ corresponds to toric Fano variety $X_P$, and -- considering polytopes up to the action of $\GL(N)$, and toric varieties up to isomorphism -- this correspondence is bijective. Our focus in this paper will be on Fano polytopes, however many of the combinatorial arguments generalise if we drop the requirement that the vertices are primitive; such polytopes correspond to toric Deligne--Mumford stacks (see, for example,~\cite{Tyo12}). For an overview of Fano polytopes, see~\cite{KN12}. Assumption~\eqref{item:generate_lattice}, although essential when considering qG-deformations of the toric variety $X_P$, is not part of the definition of Fano polytope. In two dimensions every Fano polygon satisfies assumption~\eqref{item:generate_lattice}, since every Fano polygon contains a basis for the lattice; this is not true in higher dimensions.

\begin{Ex}
Consider the three-dimensional Fano polytopes
\begin{align*}
P:=&\sconv{(1,0,0),(0,1,0),(-1,-1,0),(1,2,1),(-1,-2,-1)},\\
Q:=&\sconv{(1,0,0),(0,1,0),(-1,-1,0),(1,2,3),(-1,-2,-3)}.
\end{align*}
We have that $\abs{P\cap N}=\abs{Q\cap N}=6$, with the only non-vertex lattice point in each case corresponding to the origin. The points of $Q$ generate the sublattice $e_1\cdot\ZZ+e_2\cdot\ZZ+3e_3\cdot\ZZ$ of index three, and restricting $Q$ to this sublattice gives $P$. The Laurent polynomials
\[
f:=x+y+\frac{1}{xy}+xy^2z+\frac{1}{xy^2z},\qquad g:=x+y+\frac{1}{xy}+xy^2z^3+\frac{1}{xy^2z^3},
\]
which have Newton polytopes $P$ and $Q$ respectively, generate the same period sequence:
\begin{align*}
\pi_f(t)=\pi_g(t)&=\sum_{k=0}^\infty \sum_{m=0}^\infty {2k+3m\choose k,k,m,m,m}t^{2k+3m} \\
&=1+2t^2+6t^3+6t^4+120t^5+110t^6+1260t^7+\ldots
\end{align*}
This agrees with the period sequence for $\PP^1\times\PP^2$ given by Givental~\cite{Giv98}, and in fact the toric variety $X_P$ is equal to $\PP^1\times\PP^2$. The anti-canonical degree of $X_P$ is $(-K_{X_P})^3=\Vol{P^*}=54$, where
\[
P^*:=\{u\in\MQ\mid \text{$u(v)\geq -1$ for all $v\in P$} \}.
\]
is the \emph{dual} (or \emph{polar}) polytope to $P$ in the dual lattice $M:=\Hom(N,\ZZ)$, and $\Vol{\,\cdot\,}$ denotes the lattice-normalised volume. The degree of $X_Q$ differs from that of $X_P$ by a factor of three, that being the index of the sublattice, so that that $(-K_{X_Q})^3=\Vol{Q^*}=1/3\Vol{P^*}=18$. In particular, even though $\pi_f=\pi_g=\widehat{G}_{\PP^1\times\PP^2}$, $X_Q$ cannot be a qG-degeneration of $\PP^1\times\PP^2$.
\end{Ex}

\subsection{Vertex ansatz}
Let $P\subset\NQ$ be a \emph{smooth} Fano polytope. That is, for each facet $F$ of $P$, the vertices of $F$ generate the lattice $N$. In two dimensions there are $5$ smooth Fano polygons; in three dimensions there are only $18$ smooth Fano polytopes~\cite{Bat81,WW82}. Since $P$ is a smooth Fano polytope, the corresponding toric variety $X_P$ is a Fano manifold. Any toric Fano variety can be expressed as a GIT quotient of the form $\CC^r \GIT (\Cstar)^k$ -- this generalises the construction of projective space $\PP^{r-1}$ as $(\CC^r \setminus \{0\})/\Cstar = \CC^r \GIT \Cstar$, and amounts to specifying\footnote{In general to specify a GIT quotient $V \GIT G$, where $V$ is a vector space, we would also need to specify a character $\chi$ of $G$ called the stability condition. In our situation there is a canonical choice for $\chi$, given by $B_1+\cdots+B_r$.} the characters $B_1,\ldots,B_r$ of $(\Cstar)^k$ that define the action of $(\Cstar)^k$ on $\CC^r$. The quantum period $\widehat{G}_{X_P}$ can be computed directly from the GIT data~\cite{Giv98}. It was shown by Batyrev~\cite{Bat04} and by Batyrev--Ciocan-Fontanine--Kim--van Straten~\cite{BCKvS98} that the Laurent polynomial
\[
f_P:=\sum_{v\in\V{P}}\bx^v
\]
given by assigning the coefficient $1$ to the vertices of $P$ is a mirror partner to $X_P$, with $\pi_{f_P}=\widehat{G}_{X_P}$; this ansatz for $f_P$ is explicitly described by Przyjalkowski~\cite{Prz07}.

\begin{Prop}Let $P$ be a smooth Fano polytope. Then $f_P$ as given above is a mirror partner to $X_P$.
\end{Prop}

\begin{proof}
The quantum period of $X_P$ is computed in the following way. Given a GIT presentation of $X_P$ as above, we can write the characters $B_1,\ldots,B_r$ as a $k \times r$ integer matrix 
\[
    B=\begin{pmatrix}B_1 & \cdots & B_r\end{pmatrix}
\] 
called the weight matrix. Here $B_i$ is the weight of a torus-invariant divisor $D_i \subset X_P$. We have that
\[
\widehat{G}_{X_P}(t)=\!\!\sum_{v:\langle v,B_i\rangle\ge 0\forall i}\!\!\frac{(\sum_i\langle v,B_i\rangle)!}{\prod_i\langle v,B_i\rangle!}t^{\sum_i\langle v,B_i\rangle}.
\]
On the other hand, the classical period of $f_P$ can be expressed in terms of the matrix $A$ whose columns are the vertices of $P$ (corresponding to the torus-invariant divisors $D_i$):
\[
\pi_f(t)=\!\!\sum_{u\in \ker{A}:u_i\ge 0\forall i}\!\!\frac{(\sum_i u_i)!}{\prod_i u_i!}t^{\sum_i u_i}.
\]
We see that the sums coincide if we can identify summands via $u=vB$, but this follows directly from the facts that $A$ and $B^T$ are Gale dual matrices and that $u$ is in $\ker{A}$.
\end{proof}

\subsection{Binomial ansatz}
Consider the case when the Fano polytope $P\subset\NQ$ is \emph{reflexive}, so that $P^*\subset\MQ$ is also a Fano polytope, but in the lattice $M$. Every smooth Fano polytope is reflexive. In two dimensions there are exactly $16$ reflexive polygons, and it was observed by Galkin~\cite{Gal08} and Przyjalkowski~\cite{Prz09} that the Laurent polynomial $f_P$ obtained by assigning binomial coefficients to the monomials represented by the lattice points along the edges of $P$ is always a mirror partner to a nonsingular del~Pezzo surface. Fix an orientation on $P$, label the vertices $a_1,\ldots,a_m$ of $P$ in cyclic order, starting at some arbitrarily chosen vertex $a_1$, and set $a_{m+1}:=a_1$. For an edge $E_i=\sconv{a_i,a_{i+1}}$ of $P$, let $k_i:=\abs{E_i\cap N}-1$ denote the lattice length. Define
\begin{equation}\label{eq:binomial_edge_coeffs}
f_P:=\sum_{i=1}^m\bx^{a_i}\left(1+\bx^{w_i}\right)^{k_i}-\sum_{i=1}^m\bx^{a_i},\qquad\text{where $w_i:=\frac{1}{k_i}(a_{i+1}-a_i)\in N$ is primitive.}
\end{equation}
The second sum here removes duplicate contributions from the vertices. Explicit computation gives the following:

\begin{Prop}\label{prop:binomial_edge_coeffs}
Let $P\subset\NQ$ be a reflexive polygon. The Laurent polynomial $f_P$ defined by~\eqref{eq:binomial_edge_coeffs} is a mirror partner to a nonsingular del~Pezzo surface. Moreover, $8$ of the $10$ nonsingular del~Pezzo surfaces have a mirror partner arising in this way.
\end{Prop}

\noindent The two missing del~Pezzo surfaces here are those of the lowest degree, one and two.

The binomial ansatz is less successful in three dimensions. When applied to the three-dimensional reflexive polytopes it generates mirror partners to 92 of the 105 three-dimensional Fano manifolds, but it also generates more than 2000 Laurent polynomials which are not mirror to any three-dimensional Fano manifold. A more successful recipe in three dimensions is the Minkowski ansatz, which we now describe.

\subsection{Minkowski ansatz}
Let $P\subset\NQ$ be a three-dimensional reflexive polytope. There are $4319$ such polytopes~\cite{KS98}, and a surprisingly effective partial answer to Question~\ref{qu:assign_laurent} is given by the \emph{Minkowski ansatz}~\cite[Section~6]{ProcECM}. A toric singular point of $X_P$ corresponds to a facet $F$ of $P$, and Altmann~\cite{Alt97} tells us that deformation components of that singularity correspond to Minkowski decompositions of $F$. Altmann's work motivates the following definition. Each facet $F$ of $P$ can be decomposed into irreducible Minkowski summands:
\[
F=Q_1+\cdots+Q_r,\qquad\text{where $\dim{Q_i}\geq 1$.}
\]
We require that each $Q_i$ is either a line segment of lattice length one, or is $\GL(2,\ZZ) \rtimes \ZZ^2$-equivalent to a triangle $A_n:=\sconv{\orig,e_1,n\cdot e_2}$. If this is the case, we call the decomposition \emph{admissible}. Laurent polynomials are assigned to the summands in the obvious way: if $Q_i=\sconv{a_0,a_0+a_1}$, where $a_0,a_1\in N$, $a_1$ primitive, then $f_{Q_i}=\bx^{a_0}(1+\bx^{a_1})$; if $Q_i=\sconv{a_0,a_0+a_1,a_0+n\cdot a_2}$, where $a_0,a_1,a_2\in N$, $a_1$ and $a_2$ primitive, then $f_{Q_i}=\bx^{a_0}\left(\bx^{a_1}+(1+\bx^{a_2})^n\right)$. Define
\[
f_F:=\prod_{i=1}^rf_{Q_i}.
\]
If every facet of $P$ admits an admissible decomposition, the Laurent polynomial $f_P$ is given by the ``union'' of the $f_F$:
\[
f_P:=\sum_{v\in\bdry{P}\cap N}c_v\bx^v,\qquad\text{where $c_v$ is the coefficient of $\bx^v$ in one of the $f_F$.}
\]
Note that the definition of the $f_{Q_i}$ guarantees that the coefficients assigned to the codimension $>1$ faces agree, with binomial coefficients along the edges, so this construction is well-defined. Note also that the $f_F$ (and hence $f_P$) depend on the particular choice of admissible decomposition: different decompositions of $F$ can assign different coefficients to the monomials corresponding to $\intr{F}\cap N$.

\begin{Ex}
Consider the three-dimensional reflexive polytope with seven vertices given by
\[
P:=\sconv{(0,1,-1),(1,1,-1),(1,0,-1),(0,-1,-1),(-1,-1,-1),(-1,0,-1),(0,0,1)}\subset\NQ.
\]
The facets consist of a hexagon $F$ and six $A_1$-triangles. 
\begin{figure}[ht]
    \begin{center}
        \tdplotsetmaincoords{78}{105}
        \begin{tikzpicture}[line join=bevel,tdplot_main_coords,scale=1.2]
            \draw[black, semithick] (0,0,1) -- (0,1,-1) -- (1,1,-1) -- cycle;
            \draw[black, semithick] (0,0,1) -- (1,1,-1) -- (1,0,-1) -- cycle;
            \draw[black, semithick] (0,0,1) -- (1,0,-1) -- (0,-1,-1) -- cycle;
            \draw[black, ultra thin] (0,0,1) -- (0,-1,-1) -- (-1,-1,-1) -- cycle;
            \draw[black, ultra thin] (0,0,1) -- (-1,-1,-1) -- (-1,0,-1) -- cycle;
            \draw[black, ultra thin] (0,0,1) -- (-1,0,-1) -- (0,1,-1) -- cycle;
            
        \end{tikzpicture}
    \end{center}
\end{figure}

\noindent There are two admissible decompositions of $F$:
\begin{enumerate}
\item\label{item:hexagon_decomp_1}
the Minkowski sum of three line segments of lattice length one;
\begin{center}
\begin{tikzpicture}[scale=0.5]
\draw (-0.2,0.5) -- (0.8,0.5);
\draw (4,0) -- (5,1);
\draw (2.5,0) -- (2.5,1);
\draw (6,-0.5) -- (7,-0.5) -- (8,0.5) -- (8,1.5) -- (7,1.5) -- (6,0.5) -- (6,-0.5);

\node (1) at (1.5,0.5){$+$};
\node (2) at (3.5,0.5){$+$};
\node (3) at (5.5,0.5){$=$};
\end{tikzpicture}
\end{center}
\item\label{item:hexagon_decomp_2}
the Minkowski sum of two $A_1$-triangles.
\begin{center}
\begin{tikzpicture}[scale=0.5]
\draw (0,0) -- (1,0) -- (1,1) -- (0,0);
\draw (2,0) -- (3,1) -- (2,1) -- (2,0);
\draw (4,-0.5) -- (5,-0.5) -- (6,0.5) -- (6,1.5) -- (5,1.5) -- (4,0.5) -- (4,-0.5);

\node (1) at (1.5,0.5){$+$};
\node (1) at (3.5,0.5){$=$};
\end{tikzpicture}
\end{center}
\end{enumerate}
The Minkowski ansatz gives the Laurent polynomial
\[
f_c:=\frac{y}{z}+\frac{xy}{z}+\frac{x}{z}+\frac{1}{yz}+\frac{1}{xyz}+\frac{1}{xz}+\frac{c}{z}+z,
\]
where $c=2$ in case~\eqref{item:hexagon_decomp_1}, and $c=3$ in case~\eqref{item:hexagon_decomp_2}. The two period sequences are:
\[
\pi_{f_2}(t)=1+4t^2+60t^4+1120t^6+\ldots\qquad\text{ and }\qquad\pi_{f_3}(t)=1+6t^2+90t^4+1860t^6+\ldots
\]
The corresponding Fano manifolds are the hypersurface $\MM{2}{32}$ of bidegree~$(1,1)$ in $\PP^2\times\PP^2$, and $\MM{3}{27} = \PP^1\times\PP^1\times\PP^1$. Both Fano manifolds have a qG-degeneration to~$X_P$. The notation $\MM{k}{n}$ here indicates the $n$th three-dimensional Fano manifold of Picard rank $k$ as classified by Mori--Mukai~\cite{Mori--Mukai}, with the ordering as in~\cite{QC105}.
\end{Ex}

\subsection{Mutation of Laurent polynomials}\label{subsec:mut_Laurent}
Let $f\in\CC[\bx^{\pm1},y^{\pm1}]$, $a\in\CC[\bx^{\pm1}]$, and define the map
\[
\mu:\CC(\bx,y)\rightarrow\CC(\bx,y),\qquad\bx^vy^n\mapsto\bx^v a^n y^n.
\]
If $g:=\mu(f)$ then an application of the change-of-variables formula to the period integral~\eqref{eq:period_integral} gives
\[
\pi_f=\pi_g.
\]
Although $g$ is a rational function, in general it need not be a Laurent polynomial. If $g$ is a Laurent polynomial and if $f$ is a mirror partner to $X$, then $g$ is also a mirror partner for $X$. Write
\[
f=\sum_{i\in\ZZ}P_i(\bx)y^i\qquad\text{for some $P_i\in\CC[\bx^{\pm1}]$,}
\]
where all but finitely many of the $P_i$ are zero. Then $g$ is a Laurent polynomial if and only if for each $i\in\ZZ_{<0}$, there exists $r_i\in\CC[\bx^{\pm1}]$ such that $P_i=r_i a^{|i|}$; in this case
\[
g=\mu(f)=\sum_{i\in\ZZ_{<0}}r_i y^i+\sum_{j\in\ZZ_{\geq 0}} P_j a^j y^j\in\CC[\bx^{\pm1},y^{\pm1}].
\]
With this in mind, we make the following definition.
\begin{Defn}
Let $N$ be a lattice and let $w\in M$ be a primitive vector in the dual lattice. Then $w$ induces a grading on $\CC[N]$. Let $a\in\CC[w^\perp\cap N]$ be a Laurent polynomial in the zeroth piece of $\CC[N]$, where $w^\perp\cap N=\{v\in N\mid w(v)=0\}$. The pair $(w,a)$ defines an automorphism of $\CC(N)$ given by
\[
\mu_{w,a}:\CC(N)\rightarrow\CC(N),\qquad\bx^v\mapsto\bx^v a^{w(v)}.
\]
Let $f\in\CC[N]$. We say that $f$ is \emph{mutable with respect to $(w,a)$} if
\[
g:=\mu_{w,a}(f)\in\CC[N],
\]
in which case we call $g$ a \emph{mutation of $f$} and $a$ a \emph{factor}.
\end{Defn}

If $a=\bx^v\in\CC[w^\perp\cap N]$ is a monomial then $\mu_{w,a}$ is simply a monomial change of basis, so we regard the mutation as trivial. Similarly, we regard $(w,a)$ and $(w,\bx^v a)$, $\bx^v\in\CC[w^\perp\cap N]$, as defining equivalent mutations, that is, the factor $a$ is considered only up to ``translation''. Typically if two Laurent polynomials $f$ and $g$ are related via a change of basis then we are inclined to regard them as equivalent, however we shall see below that it can be important to remember how they were obtained.

\begin{Ex}\label{eg:Laurent_P2}
Following~\cite{HP10} consider the Laurent polynomial $f=x+y+1/(xy)\in\CC[x^{\pm1},y^{\pm1}]$. Write:
\[
x+y+\frac{1}{xy}=x\left(1+\frac{1}{x^2y}\right)+y=\frac{1}{xy}(1+xy^2)+x=x\left(1+\frac{y}{x}\right)+\frac{1}{xy}.
\]
There are three different mutations of $f$, given by:
\begin{center}
\vspace{0.5em}
\includegraphics[scale=0.6]{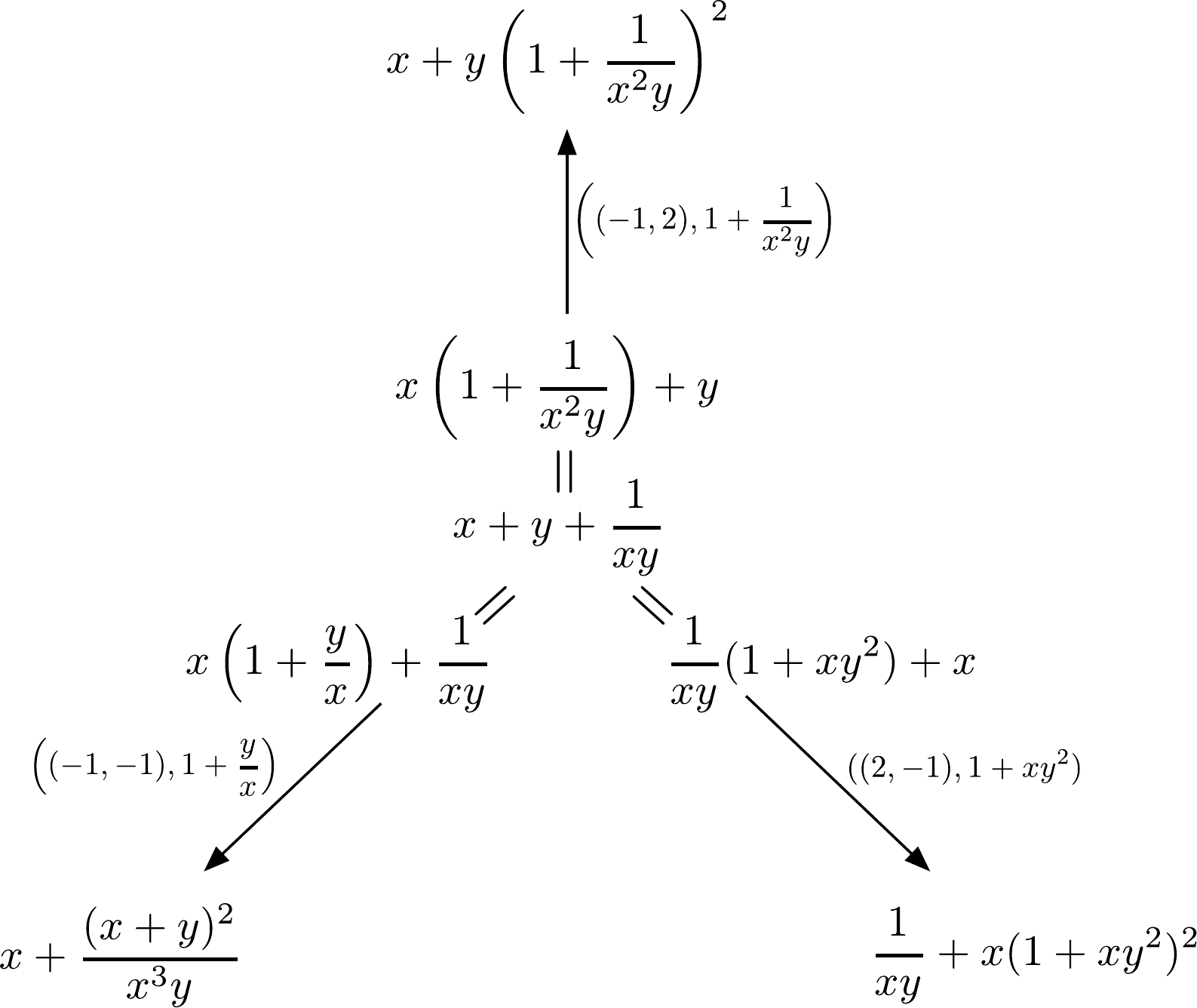}
\end{center}
Here all three mutations give equivalent results, up to monomial change of basis. However it is important to remember that we used three different mutations.
\end{Ex}

\subsection{Mutation of polytopes}\label{subsec:mut_polytope}
Let $P=\Newt{f}$ be the Newton polytope of a Laurent polynomial $f$. We want to be able to work at the level of Newton polytopes -- the ``correct'' choice of coefficients of the Laurent polynomial $f$ is addressed later. We therefore introduce the following definition, which records the effect of mutation on the Newton polytope of $f$.

\begin{Defn}
Let $N$ be a lattice and let $w\in M$ be a primitive vector in the dual lattice. Then $w$ induces a grading on $\NQ$. Let $A\subset w^\perp:=\{v\in\NQ\mid w(v)=0\}$ be a convex lattice polytope. Let $P\subset\NQ$ be a convex lattice polytope and define
\[
P_i:=\sconv{a\in P\cap N\mid w(a)=i},\qquad\text{noting that often $P_i=\emptyset$.}
\]
We say that $P$ is \emph{mutable with respect to $(w,A)$} if there exist convex lattice polytopes $R_i\subset w^\perp$ (allowing the possibility that $R_i=\emptyset$) such that
\[
\{v\in\V{P}\mid w(v)=i\}\subseteq R_i+|i| A\subseteq P_i\qquad\text{for each $i\in\ZZ_{<0}$.}
\]
If this is the case, then the \emph{mutation} of $P$ is given by the convex lattice polytope
\[
\mu_{w,A}(P):=\conv{\bigcup_{i\in\ZZ_{<0}}R_i\cup\bigcup_{j\in\ZZ_{\geq 0}}(P_j+jA)},
\]
and we call $A$ a \emph{factor}.
\end{Defn}

Although the existence of a choice of $\{R_i\}$ is required by the definition of mutation, it turns out that the resulting mutation does not depend on the choice made~\cite[Proposition~1]{ACGK12}. Hence we can safely omit the $\{R_i\}$ from the notation.

\begin{Rem}
Suppose that $g=\mu_{w,a}(f)$, for some $f$,~$g\in\CC[N]$, $a\in\CC[w^\perp\cap N]$. Then
\[
\Newt{g}=\mu_{w,\Newt{a}}\left(\Newt{f}\right).
\]
Conversely, if 
\begin{equation}
    \label{eq:polytopes_mutate}
    Q=\mu_{w,A}(P)
\end{equation}
then there exists some choice of $f$,~$g\in\CC[N]$ and $a\in\CC[w^\perp\cap N]$ with $\Newt{f}=P$, $\Newt{g}=Q$, and $\Newt{a}=A$, such that $g=\mu_{w,a}(f)$. Note however that for fixed $f$,~$g\in\CC[N]$ with $\Newt{f}=P$ and $\Newt{g}=Q$, the combinatorial statement \eqref{eq:polytopes_mutate} does not imply the existence of $a\in\CC[w^\perp\cap N]$ such that $g=\mu_{w,a}(f)$.
\end{Rem}

As in the case of Laurent polynomials, if $v\in w^\perp\cap N$ then $P\cong\mu_{w,v}(P)$, and more generally $\mu_{w,A+v}(P)\cong\mu_{w,A}(P)$. Thus the factor $A$ is considered only up to translations in $w^\perp \cap N$. Although we typically consider a polytope as being defined up to $\GL(N)$-equivalence, we will see below that it is useful to distinguish between different choices of mutation data.

Given the mirror symmetry perspective that we adopt, a key property of mutation is that it sends Fano polytopes to Fano polytopes:

\begin{Prop}[\!{\cite[Proposition~2]{ACGK12}}]
Let $Q=\mu_{w,A}(P)$ be a mutation of a convex lattice polytope $P\subset\NQ$. Then $P$ is a Fano polytope if and only if $Q$ is a Fano polytope.
\end{Prop}


\section{The Mutation Graph and Rigid Laurent Polynomials}

In this section we will introduce \emph{rigid maximally mutable Laurent polynomials}. This is a class of Laurent polynomials that provides mirror partners to all three-dimensional Fano manifolds, and more generally, we believe, to $\QQ$-factorial terminal Fano varieties in every dimension. Roughly speaking, rigid Laurent polynomials are those that are uniquely determined by the mutations that they admit, and maximally mutable Laurent polynomials are those that admit as many mutations as possible. To make this precise, we first introduce the \emph{mutation graph} of a Laurent polynomial $f$.

Identifying mutations of $f$ that differ by $\GL(N)$-equivalence can lead to problems if $f$ admits nontrivial automorphisms -- cf.~Example~\ref{eg:Laurent_P2}. But no non-trivial shear transformation can ever be an automorphism of $f$, and so it is natural to identify mutations that differ by a shear. 

\begin{Defn} A transformation $B\in\SL(N)$ is called a \emph{$w$-shear}, where $w \in M$, if $B|_{w^\perp} = \id$. 
\end{Defn}

Since $w$-shears act on $N$, they act on $\CC(N)$. Consider a mutation $g=\mu_{w,a}(f)$. If we multiply the factor $a$ by a monomial $\bx^b$, $b\in w^\perp\cap N$, then $\mu_{w,a\bx^b}(f)$ is related to $g$ by a $w$-shear. Thus considering $a$ up to multiplication by monomials in $\CC[w^\perp\cap N]$ gives $g$ up to the action of $w$-shears. We write
\[
\mu_{w,a \bx^{w^\perp\cap N}}(f)
\]
for the equivalence class. In particular, $g\in\mu_{w,a\bx^{w^\perp\cap N}}(f)$.

We can now define the mutation graph $\calG_f$ of $f$. This will be an undirected graph with vertices labelled by $\GL(N)$-equivalence classes of polytopes. Given a primitive vector $w \in M$ and a Laurent polynomial $a \in \CC[w^\perp \cap N]$, write 
\[
L(w,a) := \Big(\langle w\rangle,a\bx^{w^\perp\cap N}\Big).
\]
Here $\langle w\rangle$ denotes the linear span of $w$. Note that $a$ and $w$ determine the pair $L(w,a)$, but that the converse is false.

Before introducing the mutation graph, we fix normalisation conditions for the Laurent polynomials that we consider.

\begin{Defn} \label{def:normalisation}
A Laurent polynomial $f\in\CC[N]$ is \emph{normalised} if for all vertices $v$ of $\Newt{f}$, the coefficient of the monomial $\bx^v$ in $f$ is~$1$.
\end{Defn}

\begin{Convention}
From here onwards we assume that all Laurent polynomials (and all mutation factors) are normalised. Similarly, although our Laurent polynomials are defined over $\CC$, our expectation is that Laurent polynomials that are mirror to Fano manifolds have coefficients that are non-negative integers. From here onwards we will require that all Laurent polynomials (and all mutation factors) have non-negative integer coefficients. Recall our standing convention that if $f \in \CC[N]$ then the exponents of monomials in $f$ generate $N$.
\end{Convention}

\begin{Defn} \label{defn:mutation_graph}
Given a Laurent polynomial $f$, consider the graph $G$ with vertex labels that are Laurent polynomials and edge labels that are pairs $L(w,a)$, defined as follows. Write $\ell(v)$ for the label of a vertex $v \in V(G)$, and $\ell(e)$ for the label of an edge $e \in E(G)$.
\begin{enumerate}
\item
Begin with a vertex labelled by the Laurent polynomial $f$.
\item
Given a vertex $v$, set $g:=\ell(v)$. For each $(w,a)$, $\deg a>0$, such that $g$ is mutable with respect to $(w,a)$ and either:
\begin{enumerate}
\item
there does not exist an edge with endpoint $v$ and label $L(w,a)$; or
\item
for every edge $e=\overline{vv'}$ with $\ell(e)=L(w,a)$ we have that 
\[
\ell(v')\not\in\mu_{w,a\bx^{w^\perp\cap N}}(g);
\]
\end{enumerate}
pick a representative $g'\in\mu_{w,a\bx^{w^\perp\cap N}}(g)$ and add a new vertex $v'$ and edge $\overline{vv'}$ labelled by $g'$ and $L(w,a)$, respectively.
\end{enumerate}
The \emph{mutation graph} $\calG_f$ of $f$ is defined by removing the labels from the edges of $G$ and changing the labels of the vertices from $g$ to the $\GL(N)$-equivalence class of $\Newt{g}$. 
\end{Defn}

\begin{Defn}\label{defn:MMLP}
We partially order the mutation graphs of Laurent polynomials by saying that $\calG_f\prec \calG_g$ whenever there is a label-preserving injection $\calG_f \hookrightarrow \calG_g$. A Laurent polynomial $f$ is \emph{maximally mutable} (or for short, $f$ is an \emph{MMLP}) if $\Newt{f}$ is a Fano polytope, the constant term of $f$ is zero, and $\calG_f$ is maximal with respect to $\prec$.
\end{Defn}

\begin{Defn}\label{defn:rigidMMLP}
A maximally mutable Laurent polynomial $f$ is \emph{rigid} if the following holds: for all $g$ such that the constant term of $g$ is zero and $\Newt{f}=\Newt{g}$, if $\calG_f=\calG_g$ then $f=g$.
\end{Defn}

\begin{Rem}
Mutations leave the constant term of a Laurent polynomial unchanged; furthermore, the choice of constant term of $f$ does not affect any other coefficient in a mutation of~$f$. We have chosen to fix the constant term as zero in Definition~\ref{defn:MMLP} because the regularised quantum period \eqref{eq:regularised quantum period} of a Fano manifold $X$ always satisfies $c_1 = 0$. Similarly, our choice of normalisation conditions (Definition~\ref{def:normalisation}) reflects the properties expected of Laurent polynomial mirrors to Fano manifolds.
\end{Rem}

We believe that MMLPs form the correct class of Laurent polynomials to consider when searching for mirrors to Fano varieties. In particular, as we will justify below, we expect that \emph{every mirror to a Fano manifold is a rigid MMLP}.

In dimension two, the close relationship between mutations of Fano polygons and mutations of quivers\footnote{See, for example, \cite[Section~3.3]{Minimal-polygons} or \cite[Proposition~4.6]{Akhtar}.} means that rigidity of $f$ can be detected ``locally'' in the mutation graph. That is, one can determine whether $f$ is rigid by considering only mutations of $f$ itself. Define
\[
S_f:=\{(w,a) \mid \text{$f$ is mutable with respect to $(w,a)$}\}
\]
and, given any set $S$ of pairs $(w, a)$ with $w \in M$ primitive and $a \in \CC[w^\perp \cap N]$, define
\[
L_P(S):=\left\{f\in\CC[N] \; \middle| \; \begin{minipage}{6cm}
$\Newt{f}=P$ and $f$ is mutable with respect to $(w,a)$ for all $(w,a) \in S$
\end{minipage} \right\}.
\]
Then $f$ is a rigid MMLP if and only if
\begin{equation}\label{eq:locally_rigid}
L_{\Newt{f}}(S
_f)=\{f\}.
\end{equation}
We reemphasise that all of our Laurent polynomials (and mutation factors) here are normalised. We expect~\eqref{eq:locally_rigid} to characterise rigid MMLPs in all dimensions. 


\section{Maximal Mutability in Dimension 2}
Mutation and maximally mutable Laurent polynomials are particularly well-behaved in two dimensions, as we will now explain. The results in this section complete an important part of the program described in~\cite{pragmatic}, which recasts the classification of del~Pezzo surfaces with orbifold singularities in terms of mirror symmetry.

\subsection{Singularity content}

Let us consider mutations of polytopes in dimension two. Each edge of a Fano polygon $P\subset\NQ$ determines a point, which may be singular, in the toric surface $X_P$. We first describe the effect of mutation on these singularities: see~\cite{AK14,Alt97,Ilten:1pf}

\begin{Ex}\label{eg:polygon_mut}
Consider the mutation of Fano polygons depicted in Figure~\ref{fig:first_example}, where $w=(0,1)\in M$ and $A=\sconv{(0,0),(1,0)}\subset w^\perp$. This mutation has the effect of removing a line segment of length three at height~$-3$, and introducing a line segment of length two at height~$2$.

\begin{figure}[ht]
    \[
    \begin{tikzpicture}[scale=0.6]
        \draw (-2,-3) -- (-2,1) -- (-1,2) -- (2,-3) -- cycle;
        \node at (0,0) {$\circ$};
        \node at (-2,-3) {$\bullet$};
        \node at (-2,-2) {$\bullet$};
        \node at (-2,-1) {$\bullet$};
        \node at (-2,0) {$\bullet$};
        \node at (-2,1) {$\bullet$};
        \node at (-1,-3) {$\bullet$};
        \node at (-1,-2) {$\cdot$};
        \node at (-1,-1) {$\cdot$};
        \node at (-1,0) {$\cdot$};
        \node at (-1,1) {$\cdot$};
        \node at (-1,2) {$\bullet$};
        \node at (0,-3) {$\bullet$};
        \node at (0,-2) {$\cdot$};
        \node at (0,-1) {$\cdot$};
        \node at (1,-3) {$\bullet$};
        \node at (1,-2) {$\cdot$};
        \node at (2,-3) {$\bullet$};

        \node (arrow) at (3.5,0) {$\mapsto$};
        
        \draw (-2+8,-3) -- (-2+8,1) -- (-1+8,2) -- (1+8,2) -- (-1+8,-3) -- cycle;
        \node at (0+8,0) {$\circ$};
        \node at (-2+8,-3) {$\bullet$};
        \node at (-2+8,-2) {$\bullet$};
        \node at (-2+8,-1) {$\bullet$};
        \node at (-2+8,0) {$\bullet$};
        \node at (-2+8,1) {$\bullet$};
        \node at (-1+8,-3) {$\bullet$};
        \node at (-1+8,-2) {$\cdot$};
        \node at (-1+8,-1) {$\cdot$};
        \node at (-1+8,0) {$\cdot$};
        \node at (-1+8,1) {$\cdot$};
        \node at (-1+8,2) {$\bullet$};
        \node at (0+8,1) {$\cdot$};
        \node at (0+8,2) {$\bullet$};
        \node at (1+8,2) {$\bullet$};
    
    \end{tikzpicture}
    \]
    \caption{A mutation of Fano polygons.}
    \label{fig:first_example}
\end{figure}
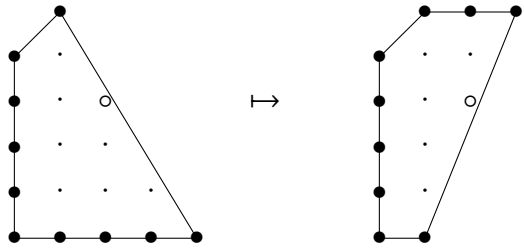

\end{Ex}

An edge $e$ of a Fano polygon $P$ determines the cone $\sigma$ over $e$, and hence a torus-fixed point in the toric variety $X_P$. A neighbourhood of this point is the affine toric variety $X_\sigma$ determined by $\sigma$, that is, the cyclic quotient singularity $\frac{1}{b}(1,c-1)$ where the rays $\{u,v\}$ of $\sigma$ are such that $u$,~$v$, and $\frac{c-1}{b} u + \frac{1}{b}v$ generate $N$. Write $b = dr$ and $c = dk$ where $d = \gcd{b,c}$. Then $X_\sigma$ is the cyclic quotient singularity $\frac{1}{dr}(1,kd-1)$. Furthermore, $r$ is the lattice height of $e$ above the origin, and $d$ is the lattice length of $e$.

\begin{figure}[ht!]
  \begin{tikzpicture}[auto]
    \draw (0,0) -- (-3*1.1,2*1.1);
    \draw (0,0) -- (3*1.1,2*1.1);
    \node at (-3,2) {$\bullet$};
    \node at (-2,2) {$\cdot$};
    \node at (-1,2) {$\bullet$};
    \node at (0,2) {$\cdot$};
    \node at (1,2) {$\bullet$};
    \node at (2,2) {$\cdot$};
    \node at (3,2) {$\bullet$};
    \node at (-1,1) {$\cdot$};
    \node at (0,1) {$\cdot$};
    \node at (1,1) {$\cdot$};
    
    \node at (3,2) {$\cdot$};
    \draw[dashed] (0,0) to (-1*1.2,2*1.2);
    \draw[dashed] (0,0) to (1*1.2,2*1.2);

    \draw[<->] (4,0) to (4,2);
    \node at (4.35,1) {$\scriptstyle r=2$};

    \draw[<->] (-3,2.7) to (3,2.7);
    \node at (0,3) {$\scriptstyle d=6$};
  \end{tikzpicture}
  \caption{A decomposition of the $T$-cone $\frac{1}{3 \cdot 2^2}(1, 5)$ into three primitive $T$-cones $\frac{1}{2^2}(1, 1)$.}
  \label{fig:T_decomp}
\end{figure}
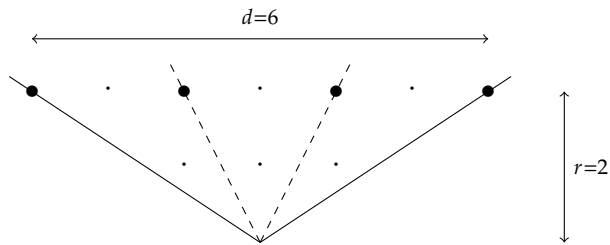

When $d=nr$ for some $n\in\ZZ$, the singularity $X_\sigma$ is called a \emph{T-singularity} and $\sigma$ is called a \emph{$T$-cone}. If $n=1$ then $X_\sigma$ is a \emph{primitive} $T$-singularity and $\sigma$ is a \emph{primitive} $T$-cone. $T$-singularities are qG-smoothable~\cite{KSB}. A general $T$-singularity $\frac{1}{n r^2}(1, knr-1)$ admits a crepant partial resolution into $n$ primitive $T$-singularities $\frac{1}{r^2}(1, kr-1)$; this amounts to the fact that, since $d = n r$, the cone $\sigma$ can be decomposed as a union of $n$ primitive $T$-cones -- see Figure~\ref{fig:T_decomp}.

\begin{figure}[b!]
  \begin{tikzpicture}[auto]
    \draw (0,0) -- (-5*1.1,3*1.1);
    \draw (0,0) -- (2*1.13,3*1.13);
    \node at (-5,3) {$\bullet$};
    \node at (-4,3) {$\cdot$};
    \node at (-3,3) {$\cdot$};
    \node at (-2,3) {$\bullet$};
    \node at (-1,3) {$\cdot$};
    \node at (0,3) {$\cdot$};
    \node at (1,3) {$\bullet$};
    \node at (2,3) {$\bullet$};
    \node at (-3,2) {$\cdot$};
    \node at (-2,2) {$\cdot$};
    \node at (-1,2) {$\cdot$};
    \node at (0,2) {$\cdot$};
    \node at (1,2) {$\cdot$};
    \node at (-1,1) {$\cdot$};
    \node at (0,1) {$\cdot$};

    \draw[dashed] (0,0) to (-2*1.13,3*1.13);
    \draw[dashed] (0,0) to (1*1.13,3*1.13);

    \draw[<->] (3,0) to (3,3);
    \node at (3.35,1.5) {$\scriptstyle r=3$};

    \draw[<->] (-5,3.7) to (2,3.7);
    \node at (-1.5,4) {$\scriptstyle d=7$};
  \end{tikzpicture}
  \caption{A decomposition of the cone $\frac{1}{7 \cdot 3}(1, 13)$ into two primitive $T$-cones $\frac{1}{3^2}(1, 5)$ and one $R$-cone $\frac{1}{3}(1,1)$.}
  \label{fig:TR_decomp}
\end{figure}
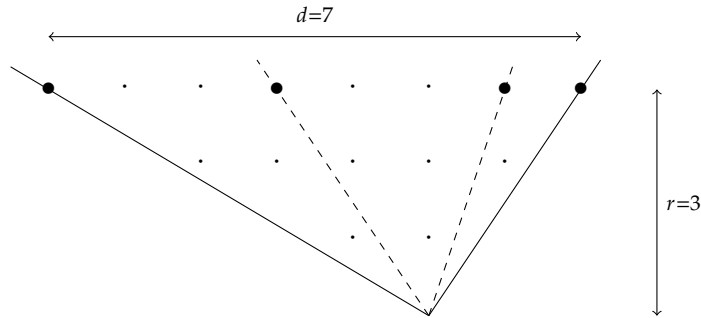

When $d = s$ with $0 < s < r$, the singularity $X_\sigma$ is called an \emph{R-singularity} and $\sigma$ is called an $R$-cone. $R$-singularities are qG-rigid. When $d = nr+s$ for some non-negative $n \in \ZZ$ and $0 < s < r$, the singularity $X_\sigma = \frac{1}{dr}(1, kd-1)$ admits a (non-unique) crepant partial resolution into $n$ primitive $T$-singularities and the single $R$-singularity $\frac{1}{sr}(1,ks-1)$; this corresponds to the (non-unique) decomposition of $\sigma$ into $n$ primitive $T$-cones and one $R$-cone -- see Figure~\ref{fig:TR_decomp}.

\noindent The singularity $X_\sigma$ qG-deforms to the R-singularity $\frac{1}{sr}(1,ks-1)$. We define the \emph{residue} of~$\sigma$ to be
\[
\res{\sigma}=\begin{cases}
\emptyset&\text{if $s=0$}\\
\frac{1}{sr}(1,ks-1)&\text{otherwise,}
\end{cases}
\]
and the \emph{singularity content} of~$\sigma$ to be the pair $(n,\res{\sigma})$.

\begin{Defn}
Let $P$ be a Fano polygon with edges $e_1,\ldots,e_m$, and suppose that the singularity content of the cone over $e_i$ is $(n_i, S_i)$. Let $n = n_1 + \cdots + n_m$, and let $\cB$ be the multiset containing those $S_i$ such that $S_i \ne \emptyset$, counted with multiplicity. The \emph{singularity content} of $P$ is
\[
\SC{P}=(n, \cB).
\]
\end{Defn}

\begin{Ex}
Returning to Figure~\ref{fig:first_example}, we see that the mutation changes the singularity content of the cone over the bottom edge from $(1, \frac{1}{3}(1,1))$ to $(0, \frac{1}{3}(1,1))$; it also introduces a new edge, at the top of the polygon, with singularity content $(1,\emptyset)$. All other singularities are unchanged. The singularity content of both polygons is $(4,\{2\times\frac{1}{3}(1,1)\})$.
\end{Ex}

\begin{Thm}[\!{\cite[Theorem~3.8]{AK14}}]
Singularity content is a mutation invariant.
\end{Thm}

If $P$ is a Fano polygon with singularity content $(n,\cB)$ then a generic qG-deformation of $X_P$ has singularities given by $\cB$. Furthermore there is a (non-unique) toric crepant partial resolution of $Y \to X_P$ such that the singularities of $Y$ are given by $n$ primitive $T$-singularities together with $\cB$.

\begin{Defn}
Let $C$ be an $R$-cone with primitive ray generators $\rho_1$,~$\rho_2 \in N$. Let $u \in M$ be the primitive normal vector to the edge $\sconv{\rho_1,\rho_2}$ that is positive on $C$, and let $r = u(\rho_1) = u(\rho_2)$. The set $\cR_C$ of \emph{residual points} of $C$ is
\[
\cR_C := C \cap \{p \in N \, \mid \, u(p) \leq r\} \setminus \{0, \rho_1, \rho_2\}.
\]
\end{Defn}

\begin{Defn}
Let $P \subset \NQ$ be a Fano polygon with singularity content $(n, \cB)$, and let $k = |\cB|$. We say that $\cR \subset N$ is \emph{a choice of residual points} for $P$ if and only if there exists a crepant subdivision of the spanning fan for $P$ into $n$ $T$-cones and $k$ $R$-cones $C_1,\ldots,C_k$ such that 
\[
\cR = \cR_{C_1} \cup \cdots \cup \cR_{C_k}.
\] 
\end{Defn}
\noindent Any two choices $\cR_1$,~$\cR_2$ of residual points for $P$ satisfy $|\cR_1| = |\cR_2|$.

\subsection{Maximally mutable Laurent polynomials} 

\begin{figure}[ht]
\[
  \begin{tikzpicture}[auto]
\fill[black!10] (0,0) -- (-1,2) -- (1,2) -- (0,0);
\draw (-2,1) -- (-1,2) -- (1,2) -- (1,-1) -- (-2,-1) -- (-2,1);
\node at (-1,2) {$\bullet$};
\node at (0,2) {$\bullet$};
\node at (1,2) {$\bullet$};
\node at (-2,1) {$\bullet$};
\node at (-1,1) {$\cdot$};
\node at (0,1) {$\cdot$};
\node at (1,1) {$\bullet$};
\node at (-2,0) {$\bullet$};
\node at (-1,0) {$\cdot$};
\node at (1,0) {$\bullet$};
\node at (-2,-1) {$\bullet$};
\node at (-1,-1) {$\bullet$};
\node at (0,-1) {$\bullet$};
\node at (1,-1) {$\bullet$};

\draw[dashed] (0,0) to (-1,2);
\draw[dashed] (0,0) to (1,2);
\draw[dashed] (0,0) to (1,-1);

\node (arrow) at (2,0) {$\mapsto$};

\fill[black!10] (0+5,0) -- (1+5,-1) -- (2+5,-1) -- (0+5,0);
\draw (-2+5,1) -- (-1+5,2) -- (2+5,-1) -- (-2+5,-1) -- (-2+5,1);
\node at (-1+5,2) {$\bullet$};
\node at (-2+5,1) {$\bullet$};
\node at (-1+5,1) {$\cdot$};
\node at (0+5,1) {$\bullet$};
\node at (-2+5,0) {$\bullet$};
\node at (-1+5,0) {$\cdot$};
\node at (1+5,0) {$\bullet$};
\node at (-2+5,-1) {$\bullet$};
\node at (-1+5,-1) {$\bullet$};
\node at (0+5,-1) {$\bullet$};
\node at (1+5,-1) {$\bullet$};
\node at (2+5,-1) {$\bullet$};

\draw[dashed] (0+5,0) to (1+5,-1);
\draw[dashed] (0+5,0) to (2+5,-1);
\draw[dashed] (0+5,0) to (-1+5,2);

  \end{tikzpicture}
\]
\caption{Mutation of a Fano polygon removes and inserts primitive $T$-cones.}\label{fig:mutation}
\end{figure}
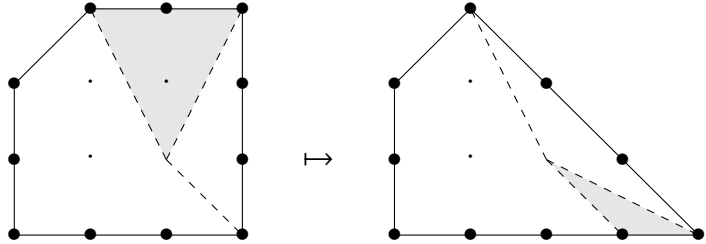

We next prove that, in dimension two, mutations of maximally mutable Laurent polynomials are in one-to-one correspondence with mutations of the underlying Newton polygons. A mutation of a polygon $P$ amounts to removing a number of primitive $T$-cones from one edge of $P$ and inserting the same number of primitive $T$-cones on the opposite side: see Figure \ref{fig:mutation}. 

Our first result produces, for each Fano polygon $P$, a non-empty family of Laurent polynomials that is compatible with all mutations of $P$. Later we will show that these Laurent polynomials are in fact maximally mutable.

\begin{Prop} \label{prop:local_MMLP}
Let $P$ be a Fano polygon. For each edge $e$ of $P$, fix a primitive vector $v_e$ in the direction of $e$ and write $w_e$ for the primitive inward-pointing normal vector to $e$. Write $(n_e, S_e)$ for the singularity content of the cone over $e$, and set $a_e = (1+\bx^{v_e})^{n_e}$. There exists a Laurent polynomial~$f$ with Newton polytope $P$ and zero constant term that is mutable with respect to $(w_e, a_e)$ for each edge~$e$. Furthermore, the free parameters in the coefficients of the general such $f$ are in bijection with any choice of residual points of~$P$, and are affine-linear functions of the coefficients of such points.
\end{Prop}

\noindent Before proving Proposition~\ref{prop:local_MMLP}, we give an example that illustrates the method of proof. 

\begin{Ex}
Consider the square $P$ with vertices $(-2,3)$, $(2,3)$, $(2,-1)$, and $(-2,-1)$. This has singularity content $(9, \{\frac{1}{3}(1,1)\})$. A normalised Laurent polynomial $f$ with Newton polygon $P$ that satisfies the mutability conditions in Proposition~\ref{prop:local_MMLP} must have binomial coefficients along the bottom edge: writing 
\[
f = \sum_{i=-1}^3 P_i(x) y^i
\]
we have that $P_{-1}$ is divisible by $(1+x)^4$ and is normalised, so $P_{-1} = x^{-2} (1+x)^4$. Furthermore, considering the top edge, we see that $P_3$ is normalised and is divisible by $(1+x)^3$, so $P_3 = x^{-2} (1+x)^4$. Arguing similarly, the other two edges also carry binomial coefficients.
\[
\begin{tikzpicture}[auto]
\draw (-2,3) -- (2,3) -- (2,-1) -- (-2,-1) -- cycle;
\node[label=above left:{\footnotesize 1}] at (-2,3) {$\bullet$};
\node[label=above:{\footnotesize 4}] at (-1,3) {$\bullet$};
\node[label=above:{\footnotesize 6}] at (0,3) {$\bullet$};
\node[label=above:{\footnotesize 4}] at (1,3) {$\bullet$};
\node[label=above right:{\footnotesize 1}] at (2,3) {$\bullet$};
\node[label=left:{\footnotesize 4}] at (-2,2) {$\bullet$};
\node[label={[label distance=-0.25cm]142.5:$\small a$}] at (-1,2) {$\cdot$};
\node[label={[label distance=-0.25cm]37.5:$\small b$}] at (0,2) {$\cdot$};
\node[label={[label distance=-0.25cm]37.5:$\small c$}] at (1,2) {$\cdot$};
\node[label=right:{\footnotesize 4}] at (2,2) {$\bullet$};
\node[label=left:{\footnotesize 6}] at (-2,1) {$\bullet$};
\node[label={[label distance=-0.25cm]142.5:$\small d$}]  at (-1,1) {$\cdot$};
\node[label={[label distance=-0.25cm]37.5:$\small e$}] at (0,1) {$\cdot$};
\node[label={[label distance=-0.35cm]57.5:$\small g$}] at (1,1) {$\cdot$};
\node[label=right:{\footnotesize 6}] at (2,1) {$\bullet$};
\node[label=left:{\footnotesize 4}] at (-2,0) {$\bullet$};
\node[label={[label distance=-0.2cm]180:$\small h$}]  at (-1,0) {$\cdot$};
\node at (0,0) {$\cdot$};
\node[label={[label distance=-0.2cm]2:\raisebox{0.1cm}{$\small k$}}] at (1,0) {$\cdot$};
\node[label=right:{\footnotesize 4}] at (2,0) {$\bullet$};
\node[label=below left:{\footnotesize 1}] at (-2,-1) {$\bullet$};
\node[label=below:{\footnotesize 4}] at (-1,-1) {$\bullet$};
\node[label=below:{\footnotesize 6}] at (0,-1) {$\bullet$};
\node[label=below:{\footnotesize 4}] at (1,-1) {$\bullet$};
\node[label=below right:{\footnotesize 1}] at (2,-1) {$\bullet$};

\draw[dashed] (0,0) to (-2,3);
\draw[dashed] (0,0) to (-1,3);
\draw[dashed] (0,0) to (2,3);
\draw[dashed] (0,0) to (2,1);
\draw[dashed] (0,0) to (2,-1);
\draw[dashed] (0,0) to (1,-1);
\draw[dashed] (0,0) to (0,-1);
\draw[dashed] (0,0) to (-1,-1);
\draw[dashed] (0,0) to (-2,-1);
\draw[dashed] (0,0) to (-2,1);

\end{tikzpicture}
\]
Now consider $P_2(x) = x^{-2}(4 + a x + b x^2 + cx^3 + 4 x^4)$. This is divisible by $(1+x)^2$, which forces
\[
4 + a x + b x^2 + cx^3 + 4 x^4 = (1+x)^2 (\alpha_1 + \alpha_2 x + \alpha_3 x^2).
\]
We must have $\alpha_1 = \alpha_3 = 4$, and $a = c = 8+\alpha$, $b = 8+2 \alpha$ for some unknown $\alpha = \alpha_2$. Similarly, $P_1(x) = x^{-2}(6 + d x + e x^2 + g x^3 + 6 x^4)$, and divisibility by $(1+x)$ forces $d = 6 + \beta$, $e = \beta + \gamma$, $g= 6 + \gamma$ for some $\beta$ and $\gamma$.

Now consider the left-hand edge, and the mutation with $w = (1,0)$. Writing
\[
f = \sum_{i=-2}^2 \bar{P}_i(y) x^i
\]
we see that $\bar{P}_{-1}(y) = y^{-1} (4+h y + d y^2 + a y^3 + 4 y^4)$ is divisible by $(1+y)^2$. This forces
\[
4 + h y + d y^2 + a y^3 + 4 y^4 = (1+y)^2 (\delta_1 + \delta_2 y + \delta_3 y^2)
\]
and therefore $\delta_1 = \delta_3 = 4$, $h = a = 8+\delta$, $d = 8 + 2 \delta$ where $\delta = \delta_2$. The relation $a = 8 + \alpha = 8 + \delta$ allows us to eliminate $\delta$, and the relation $d = 6 + \beta = 8 + 2 \delta$ allows us to eliminate $\beta$. A similar argument for the right-most edge gives us relations from $c$ and $g$; we conclude that $f$ has coefficients as follows:
\[
\begin{tikzpicture}[auto]
\draw (-2,3) -- (2,3) -- (2,-1) -- (-2,-1) -- cycle;
\node[label=above left:{\footnotesize 1}] at (-2,3) {$\bullet$};
\node[label=above:{\footnotesize 4}] at (-1,3) {$\bullet$};
\node[label=above:{\footnotesize 6}] at (0,3) {$\bullet$};
\node[label=above:{\footnotesize 4}] at (1,3) {$\bullet$};
\node[label=above right:{\footnotesize 1}] at (2,3) {$\bullet$};
\node[label=left:{\footnotesize 4}] at (-2,2) {$\bullet$};
\node[label={[label distance=-0.20cm]90:$\scriptstyle 8 + \alpha$}] at (-1,2) {$\cdot$};
\node[label={[label distance=-0.20cm]90:$\scriptstyle 8 + 2 \alpha$}] at (0,2) {$\cdot$};
\node[label={[label distance=-0.20cm]90:$\scriptstyle 8 + \alpha$}] at (1,2) {$\cdot$};
\node[label=right:{\footnotesize 4}] at (2,2) {$\bullet$};
\node[label=left:{\footnotesize 6}] at (-2,1) {$\bullet$};
\node[label={[label distance=-0.25cm]142.5:$\scriptstyle 8 + 2 \alpha$}]  at (-1,1) {$\cdot$};
\node[label={[label distance=-0.20cm]90:$\scriptstyle 4 + 4 \alpha$}] at (0,1) {$\cdot$};
\node[label={[label distance=-0.35cm]57.5:$\scriptstyle 8 + 2 \alpha$}] at (1,1) {$\cdot$};
\node[label=right:{\footnotesize 6}] at (2,1) {$\bullet$};
\node[label=left:{\footnotesize 4}] at (-2,0) {$\bullet$};
\node[label={[label distance=-0.2cm]180:\raisebox{0.1cm}{$\scriptstyle 8 + \alpha$}}]  at (-1,0) {$\cdot$};
\node at (0,0) {$\cdot$};
\node[label={[label distance=-0.2cm]2:\raisebox{0.1cm}{$\scriptstyle 8 + \alpha$}}] at (1,0) {$\cdot$};
\node[label=right:{\footnotesize 4}] at (2,0) {$\bullet$};
\node[label=below left:{\footnotesize 1}] at (-2,-1) {$\bullet$};
\node[label=below:{\footnotesize 4}] at (-1,-1) {$\bullet$};
\node[label=below:{\footnotesize 6}] at (0,-1) {$\bullet$};
\node[label=below:{\footnotesize 4}] at (1,-1) {$\bullet$};
\node[label=below right:{\footnotesize 1}] at (2,-1) {$\bullet$};

\draw[dashed] (0,0) to (-2,3);
\draw[dashed] (0,0) to (-1,3);
\draw[dashed] (0,0) to (2,3);
\draw[dashed] (0,0) to (2,1);
\draw[dashed] (0,0) to (2,-1);
\draw[dashed] (0,0) to (1,-1);
\draw[dashed] (0,0) to (0,-1);
\draw[dashed] (0,0) to (-1,-1);
\draw[dashed] (0,0) to (-2,-1);
\draw[dashed] (0,0) to (-2,1);

\end{tikzpicture}
\]
The number of free parameters here (one) is equal to the number of residual points in $P$.
\end{Ex}

\begin{figure}[bhtp]
  \begin{tikzpicture}
  
  \draw (-2,3) -- (2,3) -- (3,0);
  \draw[dotted] (-5/2,3) -- (-2,3);
  \draw[gray] (0,0) to (2,3);
  \draw[gray] (0,0) to (-1,1/2);
  \draw[gray,dotted] (-1,1/2) to (-2,1);
  
  \node at (-1.5,2+1/3) {$L$};
  \node at (-0.9,2) {$\bullet$};
  \node at (0,2) {$\bullet$};
  \node at (0.9,2) {$\bullet$};
  \node at (1.8,2) {$\bullet$};
  
  \draw[dashed] (-1.5,2) to (1.8,2);
  \draw[dotted] (4/3,2) -- (2,0);
  
  \end{tikzpicture}
  \caption{Lattice points on $L$ outside the cone $C_e$ are at $P$-height greater than $H$.}\label{fig:relative-heights}
\end{figure}
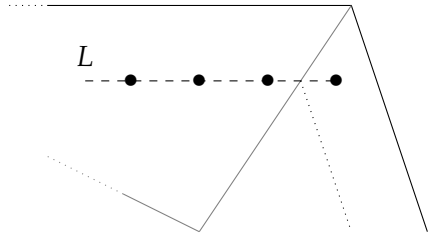

\begin{proof}[Proof of Proposition~\ref{prop:local_MMLP}]
Let us define the \emph{$P$-height} of a lattice point $p \in P$ to be the non-negative rational number $r$ such that $p$ lies on the boundary of $rP$. The set of $P$-heights of points in $P \cap N$ is a finite subset of $[0,1]$. Fix a set of residual points $\cR$ for $P$, and write $\cR_{\geq h}$ for the set of points in $\cR$ of $P$-height at least $h$. Set $f = \sum_{v \in P \cap N} a_v \bx^v$. We will prove the following statement by descending induction on $h$:
  \begin{equation}
  \tag{$\ddagger$}
  \label{induction statement}
  \begin{minipage}{0.8\textwidth}
    The coefficients $a_v$ such that $v$ has $P$-height $h$ are affine-linear functions of the coefficients $a_w$ for $w \in \cR_{\geq h}$, and these coefficients $a_w$ are independent.
  \end{minipage}
  \end{equation}
  \noindent To prove the Proposition, we need to prove statement \eqref{induction statement} for all $h \geq 0$. It holds trivially for $h=0$.

  \subsubsection*{Base case: $h=1$} Points of $P$-height $1$ lie on the boundary of $P$. Consider an edge $e$ of $P$. Without loss of generality we may assume that $e$ is horizontal and at positive vertical height $r$ above the origin. Up to overall multiplication by a monomial, the coefficients of $f$ along $e$ are
  \begin{equation}
    \label{eq:edge coefficients}
    1 + b_1 x + \cdots + b_{d-1} x^{d-1} + x^d  
  \end{equation}
  where $d = n_e r+s$ and $0 \leq s < r$, for some $b_1,\ldots,b_{d-1}$. The mutability condition gives that 
  \begin{equation}
    \label{eq:edge equation}
    1 + b_1 x + \cdots + b_{d-1} x^{d-1} + x^d = 
    (1+x)^{r n_e} \sum_{i=0}^s c_i x^i
  \end{equation}
  for some coefficients $c_i$. This forces $c_0 = c_s = 1$. If $s = 0$ or $s = 1$ then there are no residual points on $e$, and \eqref{induction statement} holds as \eqref{eq:edge coefficients} coincides with $(1+x)^d$. If $s>1$ then to prove \eqref{induction statement} it suffices to show that we can solve \eqref{eq:edge equation} uniquely for $c_i$ in terms of the coefficients $b_{k_1}, \ldots, b_{k_{s-1}}$ of residual points on $e$. That is, it suffices to show that the $(s-1) \times (s-1)$ matrix with $(i,j)$ entry
  \[
      { n_e r \choose k_i - j}
  \]
  is invertible. But the determinant of this matrix is positive -- indeed it counts certain semi-standard Young tableaux, as can be seen by specialising to $1$ all variables in the Giambelli formula~\cite[equation A6]{Fulton--Harris} that expresses Schur polynomials in terms of elementary symmetric polynomials. Thus statement \eqref{induction statement} holds for $h=1$. 

  \subsubsection*{Induction step} Suppose that $H$ satisfies $0 < H < 1$ and that statement \eqref{induction statement} holds for all $h > H$. We will prove \eqref{induction statement} for $h=H$. Let $v \in P$ be a lattice point of $P$-height $H$. Since $P$ is Fano, $v$ lies in the cone $C_e$ over a unique edge $e$. Without loss of generality we may assume that $e$ is horizontal and at positive vertical height $r$ above the origin. Let $d$ denote the lattice length of $e$, and write $d = n_e r + s$ with $0 \leq s < r$. Consider the horizontal line $L$ through $v$; this is at vertical height $Hr$ above the origin. Lattice points on $L$ fall into three classes: residual points in $C_e$, non-residual points in $C_e$, and points outside $C_e$. Convexity implies that the points on $L$ outside $C_e$ are at $P$-height greater than $H$: see Figure~\ref{fig:relative-heights}. By the induction hypothesis, therefore, coefficients of these points are fixed as affine-linear combinations of coefficients of residual points of height greater than $H$. 
  
  Up to overall multiplication by a monomial, the coefficients of $f$ along $L$ are
  \[
    b_0 + b_1 x + \cdots + b_m x^m  
  \]
  for some $m$ and $b_0,\ldots,b_m$. The mutability condition gives that 
  \begin{equation}
    \label{eq:L equation}
    b_0 + b_1 x + \cdots + b_m x^m = 
    (1+x)^{r H n_e } \sum_{i=0}^s c_i x^i
  \end{equation}
  for some coefficients $c_i$. A primitive $T$-cone with lattice height $r$ contains exactly $h$ lattice points at lattice height $h$, if $0 < h < r$, and so if $L \cap P$ contains no residual points and no points outside $C_e$ then $m<r H n_e$. In this case \eqref{eq:L equation} gives that $b_0 = \cdots = b_m = 0$, and statement \eqref{induction statement} holds. Otherwise \eqref{eq:L equation} holds with $s \geq 0$ and $s+1$ equal to the total number of points in $L \cap P$ that are either residual or lie outside $C_e$. Statement \eqref{induction statement} for $h=H$ follows if we can solve \eqref{eq:L equation} uniquely for $c_i$ in terms of the coefficients $b_{k_0}, \ldots, b_{k_s}$ of points on $L$ that are either residual or lie outside $C_e$. For this, we use the same matrix-invertibility argument as above. This completes the induction step, and the proof of Proposition~\ref{prop:local_MMLP}.
\end{proof}

\begin{Prop} \label{pro:MMLP}
  Let $F$ be Laurent polynomial in two variables with Fano Newton polygon. Then $F$ is maximally mutable if and only if it coincides with the general Laurent polynomial $f$ obtained from Proposition~\ref{prop:local_MMLP} applied to $P = \Newt{F}$.
\end{Prop}

\begin{proof}
  Suppose that $P$ is a Fano polytope, and that $f$ is a general Laurent polynomial with zero constant term and Newton polytope $P$ that satisfies the mutability conditions in Proposition~\ref{prop:local_MMLP}. These conditions imply that the mutation graph of $f$ is $n$-valent at the vertex defined by $f$, where $(n, \cB)$ is the singularity content of $P$. This is the maximum possible valency at that vertex, because $P$ admits only $n$ non-trivial mutations. Suppose now that $f$ is mutable with respect to $(w,a)$, and that the Newton polytope of the mutation $\mu_{w,a}(f)$ is~$Q$. Let $g$ be the general Laurent polynomial with zero constant term and Newton polytope~$Q$ provided by Proposition~\ref{prop:local_MMLP}. Fix residual points $\cR_P$ for $P$ and $\cR_Q$ for $Q$, noting that, because singularity content is a mutation invariant, $|\cR_P| = |\cR_Q|$. The mutation $\mu_{w,a}$ gives an invertible affine-linear identification between the coefficients of $\cR_P$ in $f$ and the coefficients of $\cR_Q$ in $g$. Thus $\mu_{w,a}f$ and $g$ coincide after a change of coefficients and, in particular, the mutation graph of $f$ is also $n$-valent at the vertex defined by $g$. Thus the mutation graph of $f$ has maximal valency at each vertex, and therefore $f$ is maximally mutable. Furthermore $\calG_F \prec \calG_f$ by construction. Since $F$ is maximally mutable it follows that $F = f$.
\end{proof}

This also proved:

\begin{Cor}  
  Let $f$ be a maximally mutable Laurent polynomial in two variables. Each vertex of the mutation graph $\calG_f$ is $n$-valent, where $(n,\cB)$ is the singularity content of $\Newt{f}$.
\end{Cor}

\begin{Ex}
  Figures~\ref{fig:mutation graph example} and~\ref{fig:quiver mutation graph example} compare the mutation graph of the rigid MMLP $f = y + \frac{(1+x)^2}{xy}$ and the quiver mutation graph~\cite{Minimal-polygons} for $\Newt{f}$. A polytope at a vertex denotes the $\GL(2,\ZZ)$-equivalence class of that polytope.

\newcommand{\polyP}[3]{
  \node (#1_0) at (#2+0,#3+0) {};
  \node (#1_1) at (#2-1/2,#3-1/2) {$\bullet$};
  \node (#1_2) at (#2+1/2,#3-1/2) {$\bullet$};
  \node (#1_3) at (#2+0,#3+1/2) {$\bullet$};
  \node (#1_mid) at (#2+0,#3-1/2) {};
  \draw (#1_1.center) -- (#1_2.center) -- (#1_3.center) -- (#1_1.center);
  \draw[dotted] (#1_0.center) -- (#1_1.center);
  \draw[dotted] (#1_0.center) -- (#1_2.center);
  \draw[dotted] (#1_0.center) -- (#1_3.center);
  \draw[dotted] (#1_0.center) -- (#1_mid.center);
}
\newcommand{\polyQ}[3]{
  \node (#1_0) at (#2+0,#3+0) {};
  \node (#1_1) at (#2-1/2,#3+0) {$\bullet$};
  \node (#1_2) at (#2+0,#3-1/2) {$\bullet$};
  \node (#1_3) at (#2+1/2,#3+0) {$\bullet$};
  \node (#1_4) at (#2+0,#3+1/2) {$\bullet$};
  \draw (#1_1.center) -- (#1_2.center) -- (#1_3.center) -- (#1_4.center) -- (#1_1.center);
  \draw[dotted] (#1_0.center) -- (#1_1.center);
  \draw[dotted] (#1_0.center) -- (#1_2.center);
  \draw[dotted] (#1_0.center) -- (#1_3.center);
  \draw[dotted] (#1_0.center) -- (#1_4.center);
}
\newcommand{\polyR}[3]{
  \node (#1_0) at (#2+0,#3+0) {};
  \node (#1_1) at (#2+1/2,#3+0) {$\bullet$};
  \node (#1_2) at (#2+0,#3+1/2) {$\bullet$};
  \node (#1_3) at (#2-9/2,#3-2/2) {$\bullet$};
  \node (#1_mid) at (#2-4/2,#3-1/2) {};
  \draw (#1_1.center) -- (#1_2.center) -- (#1_3.center) -- (#1_1.center);
  \node (#1_pt_1) at (#2-5/2,#3-1/2) {$\cdot$};
  \node (#1_pt_2) at (#2-2/2,#3+0) {$\cdot$};
  \node (#1_pt_3) at (#2-1/2,#3+0) {$\cdot$};
  \draw[dotted] (#1_0.center) -- (#1_1.center);
  \draw[dotted] (#1_0.center) -- (#1_2.center);
  \draw[dotted] (#1_0.center) -- (#1_3.center);
  \draw[dotted] (#1_0.center) -- (#1_mid.center);
}
\begin{figure}[ht!]
  \begin{tikzpicture}[scale=1.08]
    \polyP{P1}{0}{0}
    \draw (1/2,1/2) -- (5/2,5/2);
    \draw (-1/2,1/2) -- (-4/2-1/4,4/2+1/4);
    \draw (-5/2+1/4,6/2) -- (5/2-3/8,6/2);
    \draw (-3/4,-3/4) -- (-5/2,-5/2);
    \draw (3/4,-3/4) -- (6/2,-6/2);
    \polyP{P2}{-3}{3};
    \polyQ{Q}{3}{3};
    \polyR{R1}{-3}{-3};
    \polyR{R2}{6}{-3};
    \node at (0,7/2) {$\scriptstyle L\big((0,1),1+x\big)$};
    \node at (5/2,2/2+1/4) {$\scriptstyle L\big((0,1),1+x\big)$};
    \node at (-5/2-1/4,2/2+1/4) {$\scriptstyle L\big((0,1),(1+x)^2\big)$};
    \node at (5/2+1/4,-2/2-1/4) {$\scriptstyle L\big((2,-1),1+x y^2\big)$};
    \node at (-5/2-1/2,-2/2-1/4) {$\scriptstyle L\big((-2,-1),1+x y^{-2}\big)$};
    \draw[gray] (7/2+3/8, 6/2) -- (9/2+3/8, 6/2);
    \draw[gray, dotted] (9/2+3/8, 6/2) -- (10/2+3/8, 6/2);
    \draw[gray] (7/2, 7/2) -- (8/2, 8/2);
    \draw[gray, dotted] (8/2, 8/2) -- (9/2, 9/2);
    \draw[gray] (-7/2-1/4, 6/2) -- (-9/2-1/4, 6/2);
    \draw[gray, dotted] (-9/2-1/4, 6/2) -- (-10/2-1/4, 6/2);
    \draw[gray] (-7/2, 7/2) -- (-8/2, 8/2);
    \draw[gray, dotted] (-8/2, 8/2) -- (-9/2, 9/2);
    \draw[gray] (-7/2,-7/2) -- (-9/2,-9/2);
    \draw[gray, dotted] (-9/2,-9/2) -- (-10/2,-10/2);
    \draw[gray] (-5/2,-7/2) -- (-3/2,-9/2);
    \draw[gray, dotted] (-3/2,-9/2) -- (-2/2,-10/2);
    \draw[gray] (-6/2,-7/2) -- (-6/2,-9/2);
    \draw[gray, dotted] (-6/2,-9/2) -- (-6/2,-10/2);
    \draw[gray] (-4/2+13/2,-7/2-1/2) -- (-2/2+13/2,-8/2-1/2);
    \draw[gray, dotted] (-2/2+13/2,-8/2-1/2) -- (0+13/2,-9/2-1/2);
    \draw[gray] (-5/2+13/2,-7/2-1/2) -- (-3/2+12/2,-8/2-1/2);
    \draw[gray, dotted] (-3/2+12/2,-8/2-1/2) -- (-2/2+12/2,-9/2-1/2);
    \draw[gray] (-6/2+13/2,-7/2-1/2) -- (-6/2+13/2,-8/2-1/2);
    \draw[gray, dotted] (-6/2+13/2,-8/2-1/2) -- (-6/2+13/2,-9/2-1/2);
  \end{tikzpicture}
  \caption{A portion of the mutation graph of the rigid MMLP $f = y + \frac{(1+x)^2}{xy}$}
  \bigskip

  \label{fig:mutation graph example}

  \begin{tikzpicture}[scale=1.08]
    \polyP{P1}{0}{0}
    \draw (1/2,1/2) -- (5/2,5/2);
    \draw (-1/2,1/2) -- (-4/2-1/4,4/2+1/4);
    \draw (-3/4,-3/4) -- (-5/2,-5/2);
    \draw (3/4,-3/4) -- (6/2,-6/2);
    \polyQ{Q1}{-3}{3};
    \polyQ{Q2}{3}{3};
    \polyR{R1}{-3}{-3};
    \polyR{R2}{6}{-3};
    \draw[gray] (3, 7/2+3/8) -- (3, 8/2+3/8);
    \draw[gray, dotted] (3, 8/2+3/8) -- (3, 9/2+3/8);
    \draw[gray] (7/2+3/8, 6/2) -- (8/2+3/8, 6/2);
    \draw[gray, dotted] (8/2+3/8, 6/2) -- (9/2+3/8, 6/2);
    \draw[gray] (7/2, 7/2) -- (7/2+1/2/1.414, 7/2+1/2/1.414);
    \draw[gray, dotted] (7/2+1/2/1.414, 7/2+1/2/1.414) -- (7/2+1/1.414, 7/2+1/1.414);
    \draw[gray] (-7/2-1/4, 6/2) -- (-8/2-1/4, 6/2);
    \draw[gray, dotted] (-8/2-1/4, 6/2) -- (-9/2-1/4, 6/2);
    \draw[gray] (-7/2, 7/2) -- (-7/2-1/2/1.414,7/2+1/2/1.414);
    \draw[gray, dotted] (-7/2-1/2/1.414,7/2+1/2/1.414) -- (-7/2-1/1.414,7/2+1/1.414);
    \draw[gray] (-3, 7/2+3/8) -- (-3, 8/2+3/8);
    \draw[gray, dotted] (-3, 8/2+3/8) -- (-3, 9/2+3/8);
    \draw[gray] (-7/2, -7/2) -- (-7/2-1/2/1.414, -7/2-1/2/1.414);
    \draw[gray, dotted] (-7/2-1/2/1.414, -7/2-1/2/1.414) -- (-7/2-1/1.414, -7/2-1/1.414);
    \draw[gray] (-5/2,-7/2) -- (-5/2+1/2/1.414,-7/2-1/2/1.414);
    \draw[gray, dotted] (-5/2+1/2/1.414,-7/2-1/2/1.414) -- (-5/2+1/1.414,-7/2-1/1.414);
    \draw[gray] (-6/2,-7/2) -- (-6/2,-8/2);
    \draw[gray, dotted] (-6/2,-8/2) -- (-6/2,-9/2);
    \draw[gray] (-4/2+13/2,-7/2-1/2) -- (-2/2+13/2,-8/2-1/2);
    \draw[gray, dotted] (-2/2+13/2,-8/2-1/2) -- (0+13/2,-9/2-1/2);
    \draw[gray] (-5/2+13/2,-7/2-1/2) -- (-3/2+12/2,-8/2-1/2);
    \draw[gray, dotted] (-3/2+12/2,-8/2-1/2) -- (-2/2+12/2,-9/2-1/2);
    \draw[gray] (-6/2+13/2,-7/2-1/2) -- (-6/2+13/2,-8/2-1/2);
    \draw[gray, dotted] (-6/2+13/2,-8/2-1/2) -- (-6/2+13/2,-9/2-1/2);
  \end{tikzpicture}
  \caption{A portion of the quiver mutation graph for $\Newt{f}$}
  \label{fig:quiver mutation graph example}
\end{figure}

\end{Ex}

As discussed above, these results complete an important part of the program described in~\cite{pragmatic}. Fano polygons with singularity content $(n, \varnothing)$ for some $n$ fall into exactly~$10$ mutation-equivalence classes~\cite[Theorem~1.2]{Minimal-polygons}. Each mutation class supports exactly one mutation class of rigid MMLP, provided by 
Theorem~\ref{pro:MMLP}, and these correspond one-to-one with qG-deformation families of smooth del~Pezzo surfaces. Under this correspondence, the classical period $\pi_f$ of a rigid MMLP $f$ matches with the regularised quantum period $\widehat{G}_X$ of the del~Pezzo surface~\cite[\S G]{QC105}. Summarizing:

\begin{Thm}\label{thm:dim_2_rigid_smooth}
  Mutation-equivalence classes of rigid MMLPs in two variables correspond one-to-one with qG-deformation families of smooth del~Pezzo surfaces.
\end{Thm}

\noindent Combining~\cite{CKP19} with~\cite{CH17} proves the analogous result for del~Pezzo surfaces with $\frac{1}{3}(1,1)$ singularities: see~\cite{pragmatic}.

Tveiten~\cite{period-integrals} has studied the geometry of MMLPs in two dimensions. A Laurent polynomial $f$ in two variables determines a pencil of curves $Y_f \to \PP^1$ as follows. Both $f$ and the unit monomial determine sections of the anticanonical divisor on $Y_P$, where $P = \Newt{f}$ and $Y_P$ is the toric variety defined by the normal fan to $P$. Resolving singularities of $Y_P$ and resolving basepoints of the rational map $Y_P \dashrightarrow \PP^1$ defined by $[1:f]$ defines the pencil $Y_f \to \PP^1$.

\begin{Thm}[\!{\cite[Theorem~3.13]{period-integrals}}]\label{Thm:mutable-genus}Let $f$ be a maximally mutable Laurent polynomial in two variables. Then a general member of the family of curves $Y_f \to \PP^1$ has genus equal to one plus the number of residual points
  in~$\intr{P}$.
  \end{Thm}
  
Writing $p$ for the map $Y_f \to \PP^1$, we can consider the monodromy about $0 \in \PP^1$ of the local system $R^1 p_! \QQ$. Equivalently, this is the monodromy about zero of the Picard--Fuchs differential operator $L_f$ that annihilates the classical period $\pi_f$. 

\begin{Thm}[\!{\cite[Theorem~4.17]{period-integrals}}]\label{Thm:monodromy-sing-content}
  Let $f$ be a maximally mutable Laurent polynomial in two variables, let $P$ be its Newton polygon, and let $\pi_f$ be its classical period. The monodromy about zero of $L_f$ determines and is determined by the singularity content of $P$.
\end{Thm}

We have seen that singularity content is a mutation invariant of Fano polygons. It is not, however, a complete invariant.

\begin{Ex}
   Consider the polygons $P$ and $Q$ given respectively by the convex hull of the rays of a fan for $\PP^1 \times \PP^1$ and for the Hirzebruch surface $\FF_1$. Both $P$ and $Q$ have singularity content $(4,\varnothing)$, but they are not mutation equivalent. Proposition~\ref{pro:MMLP} determines unique MMLPs $f$ with Newton polytope $P$ and $g$ with Newton polytope $Q$, and these MMLPs have distinct classical periods:
   \begin{align*}
    \pi_f(t)&=1+4t^2+36t^4+400t^6+\cdots,\\
    \pi_g(t)&=1+2t^2+6t^3+6t^4+60t^5+110t^6+\cdots.
  \end{align*}
\end{Ex}

\noindent We conjecture that the classical period of the general MMLP given by Proposition~\ref{pro:MMLP} is a complete invariant for mutation of Fano polygons: cf.~\cite[Conjectures~A and~B]{pragmatic}.


\section{Some 3-dimensional Results}
The Minkowski ansatz~\cite{ProcECM} is extremely successful at recovering mirrors for the~$98$ three-dimensional Fano manifolds with very ample anticanonical bundle, but it has several drawbacks.
\begin{enumerate}
\item
The ansatz can only be applied to reflexive polytopes, and so cannot be used to recover the seven Fano manifolds that do not have very ample anticanonical bundle (although mirrors for these cases are known: see~\cite[Table~1]{QC105}).
\item
The ansatz produces~$67$ classical periods that are not the quantum period for any three-dimensional Fano manifold~\cite[\S7]{ProcECM}.
\item The ansatz is not closed under mutation. That is, there exist Laurent polynomials~$f$ mutation-equivalent to a Minkowski polynomial that are not themselves Minkowski polynomials. This holds even if one restricts attention to Laurent polynomials with reflexive Newton polytopes.
\end{enumerate}
Rigid MMLPs have none of these drawbacks. In this section we report on extensive computational experiments, which in particular give a computer-assisted proof of the following Theorem. A key ingredient here is an effective algorithm for computing the set of maximally mutable Laurent polynomials with given Newton polytope, which we will describe elsewhere~\cite{CKP-algorithm}.

\begin{Thm}\label{thm:dim3_rigid}
    Mutation-equivalence classes of rigid MMLPs $f$ such that $\Newt{f}$ is a three-dimensional reflexive polytope correspond one-to-one to the~$98$ three-dimensional Fano manifolds with very ample anticanonical bundle. Furthermore, each of the~$105$ three-dimensional Fano manifolds has a rigid MMLP mirror.
\end{Thm}
\begin{proof}
The rigid MMLPs supported on each of the $4319$ three-dimensional reflexive polytopes were computed using the computer algebra system \textsc{Magma}~\cite{Magma}. These include every Minkowski polynomial mirror to one of the $98$ three-dimensional Fano manifolds with very ample anticanonical bundle. The rigid MMLPs are listed in the Appendix; they fall into $98$ mutation-equivalence classes, which correspond one-to-one to the three-dimensional Fano manifolds just discussed. Under this correspondence, the classical period $\pi_f$ of a rigid MMLP~$f$ matches with the regularised quantum period $\widehat{G}_X$ of the corresponding Fano manifold~\cite{QC105}. Mirrors for the remaining seven Fano manifolds can be constructed from the descriptions in~\cite{QC105}. In each case, it was verified that the resulting Laurent polynomial is rigid. 
\end{proof}

One can try to produce mirrors supported on a reflexive polytope that are not rigid MMLPs. There are two obvious approaches.
\begin{enumerate}
\item
Start with a known Laurent polynomial mirror to a Fano manifold and repeatedly mutate, looking for new mirrors supported on a reflexive polytope.
\item
Start with the quantum period sequence for a Fano manifold and assign coefficients to the lattice points in a reflexive polytope in order to recover a Laurent polynomial with the correct period sequence.
\end{enumerate}
We systematically applied both approaches, but were unable to produce any non-rigid mirrors. 

Note that a three-dimensional analogue of Theorem~\ref{thm:dim_2_rigid_smooth} cannot hold: when one considers non-reflexive Newton polytopes there exist rigid MMLPs that are not mirror to a Fano manifold.

\begin{Ex} \label{ex:terminal 1}
Consider the Hilbert series with \href{http://grdb.co.uk/search/fano3?id_cmp=in&id=20522}{ID $20522$} in the database of possible Hilbert series of $\QQ$-factorial terminal Fano threefolds~\cite{GRDB-Fano3, ABR02}. Not all such potential Hilbert series are realised by genuine threefolds, but this one is: it arises from a complete intersection $X_{3,3}\subset\PP(1^5,2)$. The variety $X_{3,3}$ has a terminal singularity of type $\frac{1}{2}(1,1,1)$. A Laurent polynomial mirror to $X_{3,3}$, computed via~\cite{Giv98,HV}, is: 
\[
f=z(1+x+y)^3\left(1+\frac{1}{xyz}\right)^3-18.
\]
This is readily seen to be a rigid MMLP.
\end{Ex}

\begin{Ex} \label{ex:terminal 2}
Let $P\subset\NQ$ be the canonical Fano polytope with \href{http://grdb.co.uk/search/toricf3c?ID_cmp=in&ID=498784}{ID $498784$} in the database of toric canonical Fano threefolds~\cite{GRDB-toric3, K10}:
\[
P:=\sconv{(-1,1,-1),(1,1,-1),(-1,3,-1),(-1,-1,0),(1,-1,0),(-1,1,0),(0,0,1)}
\]
This supports a unique rigid MMLP given by:
\[
f=\left(\frac{y}{xz}+\frac{1}{xy}\right)(1+x+y)^2+z-2.
\]
By applying Laurent inversion~\cite{CKP19} we see that this defines a hypersurface of type $(2,2)$ in the four-dimensional Fano toric variety with weight data
\[
\begin{pmatrix}
1&0&1&-1&2&1\\
0&1&1&1&0&0
\end{pmatrix}
\]
Thus $f$ is a mirror to a Fano threefold that realises the Hilbert series with \href{http://grdb.co.uk/search/fano3?id_cmp=in&id=35296}{ID $35296$} in~\cite{GRDB-Fano3}; this variety has codimension~$10$ in $\PP(1^{12},2^2)$ and has~$2$ terminal singularities of type $\frac{1}{2}(1,1,1)$.
\end{Ex}

Examples~\ref{ex:terminal 1} and~\ref{ex:terminal 2} together illustrate two points that we will explore more fully elsewhere: that rigid MMLPs give a new tool~\cite{CHKP-terminal-overarching} for exploring the (unknown) classification of terminal $\QQ$-factorial Fano threefolds~\cite{ABR02}, and that they can help to find simple constructions of singular Fano threefolds that are conjectured to exist but which, as they have high codimension, otherwise seem complicated and hard to reach~\cite{H-constructions}. They also illustrate that, although proving that a given Laurent polynomial $f$ is maximally mutable is hard, proving that $f$ is \emph{rigid} maximally mutable is often much easier. Indeed any normalised Laurent polynomial $f$ with zero constant term and Fano Newton polytope such that the coefficients of $f$ are determined by the mutations that $f$ supports is automatically rigid maximally mutable.

\section{Remarks on Higher Dimension}

Every known mirror to a Fano manifold in dimension four or more is a rigid MMLP~\cite{CKP-ci,CGKS14,Kalashnikov,BGM-graph}. But when considering mirror symmetry for higher-dimensional Fano varieties, one cannot just restrict attention to smooth manifolds. For example, any simplicial terminal Fano polytope $P$ will support a unique normalised Laurent polynomial -- which is a rigid MMLP -- with $X_P$ a $\QQ$-factorial Fano variety with at worst terminal singularities. Such varieties are rigid under deformation~\cite{dFH12}. So even the foundational ansatz for mirror symmetry for toric varieties, the Givental/Hori--Vafa mirror~\cite{Giv98,HV}, leads us in higher dimensions to varieties with terminal singularities. This feature of Mirror Symmetry was masked in dimensions up to three, where the focus so far has been on Gorenstein Fano varieties: in dimension up to three, Gorenstein terminal varieties with quotient singularities are smooth. In dimensions four and higher terminal singularities are unavoidable~\cite{MMM88}.

Terminal singularities are, of course, a very natural class of singularities. Introduced by Reid, they are required by the Minimal Model Program~\cite{YPG,K13}, and the classification of terminal Fano varieties is a fundamental open problem in birational geometry. The following conjecture suggests that one approach to this classification is via rigid MMLPs; the ideas involved arose from many conversations with Alessio Corti and the rest of our collaborators in the Fanosearch program. Recall that a Fano variety $X$ is of \emph{class TG} if it admits\footnote{The only known Fano varieties that are not of class TG have $h^0({-K}_X) = 0$, and in that sense are close to being Calabi--Yau: see~\cite{CH17, Cuzzucoli20}.} a qG-degeneration with reduced fibers to a normal toric variety~\cite{pragmatic}.

\begin{Conj} \label{conj:mirror}
Rigid MMLPs in $n$ variables (up to mutation) are in one-to-one correspondence with pairs $(X,D)$, where $X$ is a Fano $n$-fold of class TG with terminal locally toric qG-rigid singularities and $D\in\abs{-K_X}$ is a general elephant (up to qG-deformation).
\end{Conj}


\subsection*{Acknowledgments}\label{sec:acknowledgements}
Much of this paper is based on research guided by and joint work with Alessio Corti. We thank him, Sergey Galkin, Vasily Golyshev, Liana Heuberger, Andrea Petracci and Thomas Prince for many useful conversations. Much of this paper was written during two visits of KT to AK in 2014, the first funded by a Scheme~4 Grant from the London Mathematical Society and the second funded by EPSRC Institutional Sponsorship~EP/N508822/1. Conjecture~\ref{conj:mirror} arose from conversations with Corti and Golyshev at the workshop \emph{Motivic Structures on Quantum Cohomology and Pencils of CY Motives} at the Max Planck Institute for Mathematics, Bonn in September 2014. 

The computations underlying this work were performed using the Imperial College High Performance Computing Service and the compute cluster at the Department of Mathematics, Imperial College London; we thank Andy Thomas and Matt Harvey for invaluable technical assistance. We thank John Cannon and the Computational Algebra Group at the University of Sydney for providing licenses for the computer algebra system \textsc{Magma}. AK is supported by EPSRC Fellowship~EP/N022513/1. TC and GP are funded by ERC Consolidator Grant~682603 and EPSRC Programme Grant~EP/N03189X/1.

\pagebreak
\begin{appendices}
\section{The Rigid Maximally Mutable Laurent Polynomials Supported on Three-Dimensional Reflexive Polytopes}

This appendix contains $100$ tables: 
\begin{center}
\begin{tabular}{rp{0.8\textwidth}}
Table~\ref{tab:summary_manifolds}: & the degree, Minkowski ID~\cite{ACGK12}, and number of rigid MMLP mirrors for each of the 98 Fano manifolds with very ample anticanonical bundle; \\
Table~\ref{tab:Fanos_for_reflexive}: & the Fano manifolds, if any, that admit rigid MMLP mirrors supported on a given three-dimensional reflexive polytope, arranged by reflexive~ID indexed from 1 to 4319~\cite{KS98};
\end{tabular}
\end{center}
\noindent and $98$ further tables, one for each three-dimensional Fano manifold $X$ with very ample anticanonical bundle, giving the rigid MMLPs that are mirror to $X$. These latter tables are hyperlinked from the corresponding entry in Table~\ref{tab:summary_manifolds}. Entries corresponding to rigid MMLPs that are not Minkowski polynomials are highlighted in $\NotMinkowski{\text{blue}}$.

\label{sec:appendices}
\begin{landscape}
\setlength{\LTcapwidth}{22cm}
\input{reflexive3_table.tex}
\end{landscape}
\end{appendices}
\bibliographystyle{amsalpha}
\bibliography{bibliography}
\end{document}